\documentclass[11pt, a4paper]{article}
\usepackage[cp1251]{inputenc}
\usepackage[russian]{babel}
\usepackage{amsmath} \usepackage{euscript}
\usepackage{mymatrix2}
\usepackage{mymatrix}
\usepackage{longtable}
\usepackage{mathrsfs}
\usepackage{amssymb}
\begin{document}

\newcounter{bnomer} \newcounter{snomer}
\newcounter{bsnomer}
\setcounter{bnomer}{0}
\renewcommand{\thesnomer}{\thebnomer.\arabic{snomer}}
\renewcommand{\thebsnomer}{\thebnomer.\arabic{bsnomer}}
\renewcommand{\refname}{\begin{center}\large{\textbf{References}}\end{center}}

\setcounter{MaxMatrixCols}{14}

\newcommand{\sect}[1]{%
\setcounter{snomer}{0}\setcounter{bsnomer}{0}
\refstepcounter{bnomer}
\par\bigskip\begin{center}\large{\textbf{\arabic{bnomer}. {#1}}}\end{center}}
\newcommand{\sst}{%
\refstepcounter{bsnomer}
\par\bigskip\textbf{\arabic{bnomer}.\arabic{bsnomer}. }}
\newcommand{\defi}[1]{%
\refstepcounter{snomer}
\par\medskip\textbf{Definition \arabic{bnomer}.\arabic{snomer}. }{#1}\par\medskip}
\newcommand{\theo}[2]{%
\refstepcounter{snomer}
\par\textbf{Теорема \arabic{bnomer}.\arabic{snomer}. }{#2} {\emph{#1}}\hspace{\fill}$\square$\par}
\newcommand{\mtheop}[2]{%
\refstepcounter{snomer}
\par\textbf{Theorem \arabic{bnomer}.\arabic{snomer}. }{\emph{#1}}
\par\textsc{Proof}. {#2}\hspace{\fill}$\square$\par}
\newcommand{\mcorop}[2]{%
\refstepcounter{snomer}
\par\textbf{Corollary \arabic{bnomer}.\arabic{snomer}. }{\emph{#1}}
\par\textsc{Proof}. {#2}\hspace{\fill}$\square$\par}
\newcommand{\mtheo}[1]{%
\refstepcounter{snomer}
\par\medskip\textbf{Theorem \arabic{bnomer}.\arabic{snomer}. }{\emph{#1}}\par\medskip}
\newcommand{\mprop}[1]{%
\refstepcounter{snomer}
\par\medskip\textbf{Proposition \arabic{bnomer}.\arabic{snomer}. }{\emph{#1}}\par\medskip}
\newcommand{\theobp}[2]{%
\refstepcounter{snomer}
\par\textbf{Теорема \arabic{bnomer}.\arabic{snomer}. }{#2} {\emph{#1}}\par}
\newcommand{\theop}[2]{%
\refstepcounter{snomer}
\par\textbf{Theorem \arabic{bnomer}.\arabic{snomer}. }{\emph{#1}}
\par\textsc{Proof}. {#2}\hspace{\fill}$\square$\par}
\newcommand{\theosp}[2]{%
\refstepcounter{snomer}
\par\textbf{Теорема \arabic{bnomer}.\arabic{snomer}. }{\emph{#1}}
\par\textbf{Схема доказательства}. {#2}\hspace{\fill}$\square$\par}
\newcommand{\exam}[1]{%
\refstepcounter{snomer}
\par\medskip\textbf{Example \arabic{bnomer}.\arabic{snomer}. }{#1}\par\medskip}
\newcommand{\deno}[1]{%
\refstepcounter{snomer}
\par\textbf{Definition \arabic{bnomer}.\arabic{snomer}. }{#1}\par}
\newcommand{\post}[1]{%
\refstepcounter{snomer}
\par\textbf{Предложение \arabic{bnomer}.\arabic{snomer}. }{\emph{#1}}\hspace{\fill}$\square$\par}
\newcommand{\postp}[2]{%
\refstepcounter{snomer}
\par\medskip\textbf{Proposition \arabic{bnomer}.\arabic{snomer}. }{\emph{#1}}%
\ifhmode\par\fi\textsc{Proof}. {#2}\hspace{\fill}$\square$\par\medskip}
\newcommand{\lemm}[1]{%
\refstepcounter{snomer}
\par\textbf{Lemma \arabic{bnomer}.\arabic{snomer}. }{\emph{#1}}\hspace{\fill}$\square$\par}
\newcommand{\lemmp}[2]{%
\refstepcounter{snomer}
\par\medskip\textbf{Lemma \arabic{bnomer}.\arabic{snomer}. }{\emph{#1}}
\par\textsc{Proof}. {#2}\hspace{\fill}$\square$\par\medskip}
\newcommand{\coro}[1]{%
\refstepcounter{snomer}
\par\textbf{Следствие \arabic{bnomer}.\arabic{snomer}. }{\emph{#1}}\hspace{\fill}$\square$\par}
\newcommand{\mcoro}[1]{%
\refstepcounter{snomer}
\par\textbf{Corollary \arabic{bnomer}.\arabic{snomer}. }{\emph{#1}}\par\medskip}
\newcommand{\corop}[2]{%
\refstepcounter{snomer}
\par\textbf{Corollary \arabic{bnomer}.\arabic{snomer}. }{\emph{#1}}
\par\textsc{Proof}. {#2}\hspace{\fill}$\square$\par}
\newcommand{\nota}[1]{%
\refstepcounter{snomer}
\par\medskip\textbf{Remark \arabic{bnomer}.\arabic{snomer}. }{#1}\par\medskip}
\newcommand{\propp}[2]{%
\refstepcounter{snomer}
\par\medskip\textbf{Proposition \arabic{bnomer}.\arabic{snomer}. }{\emph{#1}}
\par\textsc{Proof}. {#2}\hspace{\fill}$\square$\par\medskip}
\newcommand{\hypo}[1]{%
\refstepcounter{snomer}
\par\medskip\textbf{Conjecture \arabic{bnomer}.\arabic{snomer}. }{\emph{#1}}\par\medskip}
\newcommand{\prop}[1]{%
\refstepcounter{snomer}
\par\textbf{Proposition \arabic{bnomer}.\arabic{snomer}. }{\emph{#1}}\hspace{\fill}$\square$\par}

\newcommand{\Ind}[3]{%
\mathrm{Ind}_{#1}^{#2}{#3}}
\newcommand{\Res}[3]{%
\mathrm{Res}_{#1}^{#2}{#3}}
\newcommand{\epsi}{\epsilon}
\newcommand{\tri}{\triangleleft}
\newcommand{\Supp}[1]{%
\mathrm{Supp}(#1)}

\newcommand{\reg}{\mathrm{reg}}
\newcommand{\sreg}{\mathrm{sreg}}
\newcommand{\codim}{\mathrm{codim}\,}
\newcommand{\chara}{\mathrm{char}\,}
\newcommand{\rk}{\mathrm{rk}\,}
\newcommand{\chr}{\mathrm{ch}\,}
\newcommand{\id}{\mathrm{id}}
\newcommand{\Ad}{\mathrm{Ad}}
\newcommand{\col}{\mathrm{col}}
\newcommand{\row}{\mathrm{row}}
\newcommand{\low}{\mathrm{low}}
\newcommand{\pho}{\hphantom{\quad}\vphantom{\mid}}
\newcommand{\fho}[1]{\vphantom{\mid}\setbox0\hbox{00}\hbox to \wd0{\hss\ensuremath{#1}\hss}}
\newcommand{\wt}{\widetilde}
\newcommand{\wh}{\widehat}
\newcommand{\ad}[1]{\mathrm{ad}_{#1}}
\newcommand{\tr}{\mathrm{tr}\,}
\newcommand{\GL}{\mathrm{GL}}
\newcommand{\SL}{\mathrm{SL}}
\newcommand{\Sp}{\mathrm{Sp}}
\newcommand{\Mat}{\mathrm{Mat}}
\newcommand{\Stab}{\mathrm{Stab}}

\newcommand{\vfi}{\varphi}
\newcommand{\teta}{\vartheta}
\newcommand{\Bfi}{\Phi}
\newcommand{\Fp}{\mathbb{F}}
\newcommand{\Rp}{\mathbb{R}}
\newcommand{\Zp}{\mathbb{Z}}
\newcommand{\Cp}{\mathbb{C}}
\newcommand{\ut}{\mathfrak{u}}
\newcommand{\at}{\mathfrak{a}}
\newcommand{\nt}{\mathfrak{n}}
\newcommand{\mt}{\mathfrak{m}}
\newcommand{\htt}{\mathfrak{h}}
\newcommand{\spt}{\mathfrak{sp}}
\newcommand{\rt}{\mathfrak{r}}
\newcommand{\rad}{\mathfrak{rad}}
\newcommand{\bt}{\mathfrak{b}}
\newcommand{\gt}{\mathfrak{g}}
\newcommand{\vt}{\mathfrak{v}}
\newcommand{\pt}{\mathfrak{p}}
\newcommand{\Xt}{\mathfrak{X}}
\newcommand{\Po}{\mathcal{P}}
\newcommand{\Uo}{\EuScript{U}}
\newcommand{\Fo}{\EuScript{F}}
\newcommand{\Do}{\EuScript{D}}
\newcommand{\Eo}{\EuScript{E}}
\newcommand{\Iu}{\mathcal{I}}
\newcommand{\Mo}{\mathcal{M}}
\newcommand{\Nu}{\mathcal{N}}
\newcommand{\Ro}{\mathcal{R}}
\newcommand{\Co}{\mathcal{C}}
\newcommand{\Lo}{\mathcal{L}}
\newcommand{\Ou}{\mathcal{O}}
\newcommand{\Au}{\mathcal{A}}
\newcommand{\Vu}{\mathcal{V}}
\newcommand{\Bu}{\mathcal{B}}
\newcommand{\Sy}{\mathcal{Z}}
\newcommand{\Sb}{\mathcal{F}}
\newcommand{\Gr}{\mathcal{G}}

\author{Mikhail V. Ignatyev, Anton S. Vasyukhin\thanks{This work was partially supported by RFBR grant no. 13--01--97000--р\_поволжье\_а. The first author was partially supported by the Dynasty Foundation and by DAAD program ``Forschungsaufenthalte f$\ddot{\mathrm{u}}$r Hochschullehrer und Wissenschaftler'' ref. no. A/13/00032.}
}
\date{\small Samara State University\\
Chair of algebra and geometry\\
\texttt{mihail.ignatev@gmail.com}, \texttt{safian.malk@gmail.com}}
\title{\Large{Rook placements in $A_n$ and combinatorics of $B$-orbit closures}} \maketitle

\sect{Introduction and the first main result}

\sst Let $G$ be a complex reductive algebraic group, $T$ a maximal torus of~$G$, $B$ a Borel subgroup of~$G$ containing $T$, and $U$ the unipotent radical of $B$. Let $\Phi$ be the root system of $G$ with respect to $T$, $\Phi^+$ the set of positive roots with respect to $B$, $\Delta$ the set of fundamental roots, and $W$ the Weyl group of $\Phi$ (see \cite{Bourbaki}, \cite{Humphreys} and \cite{Humpreys2} for basic facts about algebraic groups and root systems). Denote by $\Fo=G/B$ the flag variety and by $X_w\subseteq\Fo$ the Schubert subvariety corresponding to an element $w$ of the Weyl group $W$. Let $\gt$, $\bt$, $\nt$ be the Lie algebras of the groups $G$, $B$, $U$ respectively, $\nt^*$ the dual space of $\nt$. The group $B$ acts on $\nt$ by the adjoint action; the dual action of $B$ on $\nt^*$ is called coadjoint. We will denote the result of this action by $b.\lambda$, $b\in B$, $\lambda\in\nt^*$. According to the orbit method discovered by A.A. Kirillov in 1962, orbits of the coadjoint action play a key role in representation theory of $B$ and~$U$ \cite{Kirillov1}, \cite{Kirillov2}.

Denote by $\leq_B$ the Bruhat--Chevalley order on $W$. This order plays a fundamental role in~a~multitude of contexts. For instance, it describes the incidences among Schubert subvarieties of the flag variety~$\Fo$, i.e., $X_v\subseteq X_w$ if and only if $v\leq_B w$. An interesting subposet of the Bruhat--Chevalley order is induced by the involutions, i.e., the elements of order 2 of $W$. We denote this subposet by $I(W)$. Activity around~$I(W)$ was initiated by R. Richardson and T. Springer \cite{RichardsonSpringer}, who proved that the inverse Bruhat--Chevalley order on~$I(S_{2n+1})$ encodes the incidences among the closed orbits under the action of the Borel subgroup of the special linear group on the symmetric variety $\mathrm{SL}_{2n+1}({\mathbb{C}})/\mathrm{SO}_{2n+1}({\mathbb{C}})$. (Here $S_n$ denotes the symmetric group on~$n$ letters). The poset of involutions was also studied by F.~In\-cit\-ti \cite{Incitti1}, \cite{Incitti2} from a purely combinatorial point of view. In particular, he described the covering relation of this poset.

To each involution $w\in I(W)$ one can assign the coadjoint $B$-orbit $\Omega_w\subseteq\nt^*$ be the following rule. The involution $w$ can be canonically expressed as a product of pairwise commuting reflections~\cite{Deodhar},~\cite{Springer}: $$w=\prod_{\beta\in D}s_{\beta},~D\subseteq\Phi^+.$$ (Since $s_{\beta}$'s commute, $D$ is an orthogonal subset of $\Phi^+$.) The root vectors $e_{\alpha}$, $\alpha\in\Phi^+$, form a basis of~$\nt$. Denote by $\{e_{\alpha}^*,~\alpha\in\Phi^*\}$ the dual basis of $\nt^*$. By definition, $\Omega_w$ is the orbit of the element $$f_w=\sum_{\beta\in D}e_{\beta}^*.$$

It was shown in \cite{Ignatyev1} that the Bruhat--Chevalley order on $I(S_n)$ encodes the incidences among the closures of such orbits for $G=\GL_n(\Cp)$. Namely, $\Omega_{\sigma}\subseteq\overline{\Omega}_{\tau}$ if and only if $\sigma\leq_B\tau$. These results are in some sense ``dual'' to A. Melnikov's results \cite{Melnikov1}, \cite{Melnikov2}. In~\cite{Ignatyev2}, similar results were obtained for $G=\mathrm{Sp}_{2n}(\Cp)$. Using results of A. Panov \cite{Panov}, the first author proved in \cite{Ignatyev1} and \cite{Ignatyev2} that for $\GL_n(\Cp)$ and $\mathrm{Sp}_{2n}(\Cp)$, $\dim\Omega_w$ equals $l(w)$, the length of $w$. Further, in \cite{Panov}, Panov constructed polarizations at the elements $f_w$, $w\in I(S_n)$. (Polarizations play an important role in the orbit method \cite{Kirillov2}.) In \cite{Ignatyev3}, \cite{Ignatyev4} similar results were obtained for other reductive groups. The goal of the paper is to generalize these results to the case of orbits associated with rook placements in the root system $A_n$; see the next Subsection for precise definitions.

The paper is organized as follows. In the next Subsection, we define rook placements in root systems. We also define a certain partial order on the set of rook placements in $A_n$ and describes connections of this poset with the $B$-orbit closures for $G=\GL_n(\Cp)$. Our first main result is Theorem~\ref{mtheo:order}.
In~Section~\ref{sect:pol_dim}, we formulate and prove our second main result, Theorem~\ref{mtheo:dim} (see Subsection~\ref{sst:dim_pol_formulir}). Precisely, we construct polarizations at some linear forms on orbits associated with rook placements and compute the dimension of the orbit associated with a rook placement in purely combinatorial terms. We also give an upper bound for the dimension in terms of the Weyl group. Finally, in Section~\ref{sect:cov_rel}, we describe the covering relation of the poset of rook placements. Our third main result, Theorem~\ref{mtheo:cov_rel}, is formulated in Subsection~\ref{sst:cov_rel_formulir} and proved in Subsections~\ref{sst:cov_ness},~\ref{sst:cov_suff}.

\medskip\textsc{Acknowledgements}. We thanks Professor Aleksandr Panov for useful discussions. The work was partially supported by the Russian Foundation for Basic Research, research project no. 13--01--97000--р\_поволжье\_а. RFBR is gratefully acknowledged. Mikhail Ignatyev also express his gratitude to the Dynasty Foundation for the financial support. The work was completed during the stay of Mikhail Ignatyev at Jacobs University, Bremen, supported by DAAD program ``For\-schungsaufenthalte f$\ddot{\mathrm{u}}$r Hochschullehrer und Wissenschaftler'' ref. no. A/13/00032. Mikhail Ignatyev thanks Professor Ivan Penkov for his hospitality and useful discussions, and DAAD for the financial support.

\sst Here we give precise definitions and prove our first main result.

\defi{A \emph{rook placement} is a subset $D\subseteq\Phi^+$ such that $(\alpha,\beta)\leq0$ for all $\alpha$, $\beta\in D$.}

\exam{Let $G=\GL_n(\Cp)$ be the general linear group, so $\Phi=A_{n-1}$. As usual, we identify~$\Phi^+$ with the set $$\{\epsi_j-\epsi_i,~1\leq j<i\leq n\}\subset\Rp^n,$$ where $\{\epsi_i\}_{i=1}^n$ is the standard basis. For instance, if $n=8$, then $$D=\{\epsi_1-\epsi_3,~\epsi_2-\epsi_6,~\epsi_3-\epsi_7,~\epsi_4-\epsi_5,~\epsi_6-\epsi_8\}$$ is a rook placement. Here we draw this subset on the lower-triangular chessboard using the following rule: the $(i,j)$th box is occupied by a rook if and only if $\epsi_j-\epsi_i\in D$.
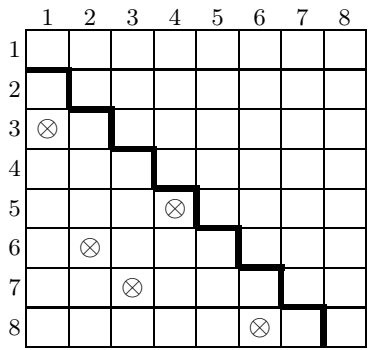
\begin{equation*}
\mymatrix{ \pho& \pho& \pho& \pho& \pho& \pho& \pho& \pho\\
\Top{2pt}\Rt{2pt}\pho& \pho& \pho& \pho& \pho& \pho& \pho& \pho\\
\otimes& \Top{2pt}\Rt{2pt}\pho& \pho& \pho& \pho& \pho& \pho& \pho\\
\pho& \pho& \Top{2pt}\Rt{2pt}\pho& \pho& \pho& \pho& \pho& \pho\\
\pho& \pho& \pho& \Top{2pt}\Rt{2pt}\otimes& \pho& \pho& \pho& \pho\\
\pho& \otimes& \pho& \pho& \Top{2pt}\Rt{2pt}\pho& \pho& \pho& \pho\\
\pho& \pho& \otimes& \pho& \pho& \Top{2pt}\Rt{2pt}\pho& \pho& \pho\\
\pho& \pho& \pho& \pho& \pho& \otimes& \Top{2pt}\Rt{2pt}\pho& \pho\\
}\end{equation*}
Note that the definition of a rook placement guarantees that rooks don't hit each other.\label{exam:RP1}}

To each rook placement $D\subseteq\Phi^+$ one can assign the $B$-orbit $\Omega_D\subseteq\nt^*$ by putting $$f_D=\sum_{\beta\in D}e_{\beta}^*$$ and $\Omega_D=B.f_D$. (Clearly, if $D$ is an orthogonal subset, then $\Omega_D=\Omega_w$ for $w=\prod_{\beta\in D}s_{\beta}$.) Further, if $\xi\colon D\to\Cp^{\times}$ is a map, then we put $\Theta_{D,\xi}=U.f_{D,\xi}$, where $$f_{D,\xi}=\sum_{\beta\in D}\xi(\beta)e_{\beta}^*.$$ We say that the orbits $\Omega_D$ and $\Theta_{D,\xi}$ are \emph{associated} with the rook placement $D$. In fact,
$$\Omega_D=\bigcup_{\xi\colon D\to\Cp^{\times}}\Theta_{D,\xi}$$
for $G=\GL_n(\Cp)$, see Lemma~\ref{lemm:Omega_union_Theta} below (cf. \cite[Theorem 2.5]{Kostant1}, \cite[Corollary 3.3]{Melnikov1} and
\cite[Lemma 2.1]{Ignatyev1}). Note that almost all coadjoint $U$-orbits studied up to now are associated with certain rook placements, see, e.g., \cite{Andre1}, \cite{Andre2}, \cite{AndreNeto1}, \cite{Kirillov3}, \cite{Kostant1}, \cite{Kostant2} and \cite{Panov}.

Form now on, let $G$ be the general linear group $\GL_n(\Cp)$, so $W\cong S_n$, the symmetric group on $n$ letters. Denote by $B$ the group of all invertible upper-triangular matrices, then $U$ is the unitriangular group, i.e., the group of all upper-triangular matrices with 1's on the diagonal, and $\nt$ is the space of all upper-triangular matrices with zeroes on the diagonal. We identify $\Phi^+=A_{n-1}^+$ with the set
$$\{(i,j)\in\Zp\times\Zp\mid 1\leq j<i\leq n\}$$
by sending $$\epsi_j-\epsi_i\mapsto(i,j).$$ Now, if $\alpha=(i,j)\in\Phi^+$ is a root, then the root vector $e_{\alpha}$ is an elementary matrix, namely, $e_{\alpha}=e_{j,i}$. Using the trace form
$$\langle\lambda,x\rangle=\tr{\lambda x},~\lambda\in\nt^t,~x\in\nt,$$
one can identify $\nt^*$ with the space $\nt^t$ of all lower-triangular matrices with zeroes on the diagonal. Under this identification, $e_{j,i}^*=e_{i,j}$, so
\begin{equation*}
(f_D)_{i,j}=\begin{cases}
1,&\text{if }(i,j)\in D,\\
0&\text{otherwise.}
\end{cases}
\end{equation*}
Note also that the coadjoint action has the simple form
\begin{equation*}
b.\lambda=(b\lambda b^{-1})_{\low},~b\in B,~\lambda\in\nt^*,
\end{equation*}
where $X_{\low}$ denotes the strictly lower-triangular part of a matrix $X\in\Mat_n(\Cp)$.

Denote by $T$ the group of invertible $n\times n$ diagonal matrices. Recall that $B$ is a semi-direct product of $U$ and $T$, $B=U\rtimes T$. In particular, for any $g\in B$, there exist $u\in U$, $t\in T$ such that $g=ut$. Denote by $1_n$ the $n\times n$ identity matrix. Finally, suppose $$D=\{(i_1,j_1),\ldots(i_s,j_s)\},$$ $i_l>j_l$, $j_1<\ldots<j_s$, is a rook placement. Then put
\begin{equation}
w=(i_1,j_1)\ldots(i_s,j_s)\in S_n.\label{formula:w_D}
\end{equation}
In other words, $w=s_{\beta_1}\ldots s_{\beta_s}$, where $\beta_l=(i_l,j_l)$. We need the following simple Lemma (cf. \cite[Lemma 2.1]{Ignatyev1}).

\newpage\lemmp{Let $D$ be a rook placement. Then $\Omega_D=\bigcup\nolimits_{\xi\colon D\to\Cp^{\times}}\Theta_{D,\xi}$.\label{lemm:Omega_union_Theta}}{Let
$\xi\colon D\to\Cp^{\times}$ be a map. Suppose $D=\{(i_1,j_1),\ldots,(i_s,j_s)\}$. Let $w$ be the permutation defined by (\ref{formula:w_D}). One can express $w$ as a product of disjoint cycles:
$$w=(a_1^1,\ldots,a_{b_1}^1)\ldots(a_1^k,\ldots,a_{b_k}^k),$$
$a_1^i<\ldots<a_{b_i}^i$ for all $1\leq i\leq k$. It is easy to see that
$$D=\bigcup_{i=1}^k\bigcup_{j=2}^{b_i}\{(a_j^i,a_{j-1}^i)\}.$$

Then define $t$ to be the diagonal matrix such that
\begin{equation*}
t_{p,p}=\begin{cases}
\prod_{z=2}^j\xi(a_z^i,a_{z-1}^i)^{-1},&\text{if }p=a_j^i\text{ for some }1\leq i\leq k,~1\leq j\leq b_i,\\
1&\text{otherwise}.
\end{cases}
\end{equation*}
For instance, if $n=8$ and $D=\{(3,1),(6,2),(7,3),(5,4),(8,6)\}$, as in Example~\ref{exam:RP1}, then
$$w=(1,3,7)(2,6,8)(4,5),$$
and so
$$t=\mathrm{diag}\left\{1,1,\dfrac{1}{\xi(3,1)},1,\dfrac{1}{\xi(5,4)},
\dfrac{1}{\xi(6,2)},\dfrac{1}{\xi(3,1)\cdot\xi(7,3)},\dfrac{1}{\xi(6,2)\cdot\xi(8,6)}\right\}.$$
One can trivially check that $t.f_{D,\xi}=f_D$, hence $\Theta_{D,\xi}\subseteq\Omega_D$.

On the other hand, let $g$ be an element of $B$. Then there exist
$u\in U$, $t\in T$ such that $g=ut$, so $g.f_D=u.f_{D,\chi}$,
where $\chi(i_r,j_r)=g_{i_r,i_r}/g_{j_r,j_r}$, $1\leq r\leq s$. Thus,
$g.f_D\in\Theta_{D,\chi}$.}

To discuss the incidences among the closures of orbits associated with rook placements, we need some more notation. By $\overline{Z}$ we denote the Zariski closure of a subset $Z\subseteq\nt^*$. Given $X$, $Y\in\Mat_n(\Zp)$, we write $X\leq Y$ if and only if $X_{i,j}\leq Y_{i,j}$ for all $i,j$. To each rook placement $D\subseteq\Phi^+$ we assign the matrix $R_D$ by the following rule:
\begin{equation*}
(R_D)_{i,j}=\begin{cases}
\rk\pi_{i,j}(f_D),&\text{if }i>j,\\
0&\text{otherwise,}
\end{cases}
\end{equation*}
where $\pi_{i,j}(X)$ denotes the lower-left triangular part of a matrix $X\in\Mat_n(\Cp)$. In other words, $(R_D)_{i,j}$, $i>j$, is just the number of rooks situated non-strictly to the South-West of the box $(i,j)$. Finally, we put $D\leq D'$ if $R_D\leq R_{D'}$.

\exam{Let $n=8$ and $D$ be as in the previous example. Then
\begin{equation*}
R_D=\mymatrix{\fho0& \fho0& \fho0& \fho0& \fho0& \fho0& \fho0& \fho0\\
\Top{2pt}\Rt{2pt}\fho1& \fho0& \fho0& \fho0& \fho0& \fho0& \fho0& \fho0\\
\fho1& \Top{2pt}\Rt{2pt}\fho2& \fho0& \fho0& \fho0& \fho0& \fho0& \fho0\\
\fho0& \fho1& \Top{2pt}\Rt{2pt}\fho2& \fho0& \fho0& \fho0& \fho0& \fho0\\
\fho0& \fho1& \fho2& \Top{2pt}\Rt{2pt}\fho3& \fho0& \fho0& \fho0& \fho0\\
\fho0& \fho1& \fho2& \fho2& \Top{2pt}\Rt{2pt}\fho2& \fho0& \fho0& \fho0\\
\fho0& \fho0& \fho1& \fho1& \fho1& \Top{2pt}\Rt{2pt}\fho2& \fho0& \fho0\\
\fho0& \fho0& \fho0& \fho0& \fho0& \fho1& \Top{2pt}\Rt{2pt}\fho1& \fho0\\
}\end{equation*}}

Our first result is as follows (cf. \cite[Theorem 1.10]{Ignatyev1}, \cite[Theorem 3.5]{Melnikov2}).

\newpage\mtheop{Let $D$, $D'$ be rook placements. If $\Omega_D\subseteq\overline{\Omega}_{D'}$, then $D\leq D'$.\label{mtheo:order}}{Using Lemma~\ref{lemm:Omega_union_Theta}, one can repeat literally the proofs of \cite[Lemma 2.2]{Ignatyev1} and \cite[Proposition 2.3]{Ignatyev1} to obtain the result. Namely, we claim that $$\rk\pi_{i,j}(\lambda)
=(R_D)_{i,j}$$ for all $\lambda\in\Omega_D$, $1\leq j<i\leq
n$. Indeed, Lemma \ref{lemm:Omega_union_Theta} shows
that it's enough to check that if $u\in U$, $f\in\nt^*$, then
$$\rk\pi_{i,j}(u.f)=\rk\pi_{i,j}(f)$$ for all $1\leq j<i\leq
n$, because
$$\rk\pi_{i,j}(f_{D,\xi})=(R_D)_{i,j}$$ for all maps $\xi\colon D\to\Cp^{\times}$.

To do this, pick an element $u\in U$. It's well-known that there exist
$\alpha_{j,i}\in\Cp$ such that
\begin{equation*}
u=\prod_{(i,j)\in\Phi^+}x_{j,i}(\alpha_{j,i}),
\end{equation*}
where $x_{j,i}(\alpha_{j,i})=1_n+\alpha_{j,i}e_{j,i}$ (the product
is taken in any fixed order). Hence we can assume
$u=x_{j,i}(\alpha)$ for some $(i, j)\in\Phi^+$, $\alpha\in\Cp$. Then
\begin{equation*}\predisplaypenalty=0
(u.f)_{r,s}=\begin{cases}f_{j,s}+\alpha f_{i, s},&\text{if
}r=j\text{ and }1\leq s<j,\\
f_{r,i}-\alpha f_{r, j},&\text{if
}s=i\text{ and }i<r\leq n,\\
f_{r,s}&\text{otherwise}.
\end{cases}
\end{equation*}

Hence if $r>j$ and $s<i$, then $\pi_{r,s}(u.f)=\pi_{r,s}(f)$. If
$r\leq j$ (and so $s<r\leq j<i$), then the $j$th row of
$\pi_{r,s}(u.f)$ is obtained from the $j$th row of
$\pi_{r,s}(f)$ by adding the $i$th row of $\pi_{r,s}(f)$
multiplied by $\alpha$. Similarly, if $s\geq i$ (and so $r>s\geq
i>j$), then the $i$th column of $\pi_{r,s}(u.f)$ is obtained from
the $i$th column of $\pi_{r,s}(f)$ by subtracting the $j$th
column of $\pi_{r,s}(f)$ multiplied by $\alpha$. In both cases,
$\rk\pi_{r,s}(u.f)=\rk\pi_{r,s}(f)$, as required.

Now, suppose $\Omega_D\subseteq\overline{\Omega}_{D'}$, but $D\nleq D'$. This means that there exists a root
$(i,j)\in\Phi^+$ such that $(R_{D'})_{i,j}<(R_D)_{i,j}$.
Denote
\begin{equation*}
Z=\{f\in\nt^*\mid\rk\pi_{r,s}(f)\leq(R_{D'})_{r,s}\text{ for
all }(r,s)\in\Phi\}.
\end{equation*}
Clearly, $Z$ is closed with respect to Zariski topology. We proved that $\Omega_{D'}\subseteq
Z$, hence $\overline{\Omega}_{D'}\subseteq Z$. But
$f_D\notin Z$, hence $\Omega_{D}\nsubseteq Z$, a
contradiction.}

\nota{i) It follows that $\Omega_D=\Omega_{D'}$ if and only if $D=D'$. Indeed, assume $\Omega_D=\Omega_{D'}$. Then $\Omega_D\subseteq\overline{\Omega}_{D'}$ and vice versa, so $R_D\leq R_{D'}$ and $R_{D'}\leq R_D$. Hence $R_D=R_{D'}$. One can easily deduce from this fact that $D=D'$.

ii) In \cite{Smirnov}, E.Yu. Smirnov studied $B$-orbits on the direct product of two Grassmanninans. They are indexed by the set of rook placements on certain Young diagrams. As one can see from\linebreak \cite[Theorem 3.10]{Smirnov}, the partial order on this set induced by the incidences among the closures of such orbits is closely related to the partial order defined above. It would be interesting to investigate any deeper relations between the situation considered by Smirnov and our situation.

iii) Note that if $D$ and $D'$ are orthogonal subsets of $\Phi^+$, then $D\leq D'$ implies $\Omega_D\subseteq\overline{\Omega}_{D'}$. But in general this is not true. For instance, let $n=4$, $D=\{(3,2),(4,3)\}$, $D'=\{(2,1),(3,2)\}$. Clearly, $D\leq D'$. On the other hand, one can easily check that $$\lambda_{4,2}\lambda_{2,1}+\lambda_{4,3}\lambda_{3,1}=0$$ for all $\lambda\in\Omega_{D'}$, so $f_D\notin\overline{\Omega}_{D'}$.}

There is a corollary of Theorem~\ref{mtheo:order} in terms of the Bruhat--Chevalley order on the set of so-called \emph{Kerov involutions} \cite{Kerov}. To each rook placement $D\subseteq\Phi^+$ one can assign the involution $\sigma_D\in S_{2n-2}$ by the following rule: if
$$D=\{(i_1,j_1),\ldots,(i_s,j_s)\},$$
then, by definition, $$\sigma_D=(2i_1-2,2j_1-1)\ldots(2i_r-2,2j_r-1).$$

\exam{Let $n=8$ and $D=\{(3,1),(6,2),(7,3),(5,4),(8,6)\}$, as in Example~\ref{exam:RP1}. Then
\begin{equation*}
\begin{split}
\sigma_D&=(4,1)\cdot(10,3)\cdot(12,5)\cdot(8,7)\cdot(14,11)\\
&=\begin{pmatrix}
1&2&3&4&5&6&7&8&9&10&11&12&13&14\\
4&2&10&1&12&6&8&7&9&3&14&5&13&11
\end{pmatrix}\in S_{14}.
\end{split}
\end{equation*}}

\mcorop{Let $D$, $D'$ be rook placements. If $\Omega_D\subseteq\overline{\Omega}_{D'}$, then $\sigma_D\leq_B\sigma_{D'}$.}{Denote by $\wt D$, $\wt D'$ the orthogonal subsets of $A_{2n-3}^+$ corresponding to the involutions $\sigma_D$, $\sigma_{D'}$ respectively. It follows from \cite[Theorem 1.10]{Ignatyev1} that $\sigma_D\leq_B\sigma_{D'}$ is equivalent to $R_{\wt D}\leq R_{\wt D'}$. But the last inequality is equivalent to $R_D\leq R_{D'}$, i.e., $D\leq D'$.}

\sect{Polarizations and dimensions of orbits}\label{sect:pol_dim}

\sst\label{sst:dim_pol_formulir} Our next goal is to compute the dimension of $\Omega_D$. To do this, we need some more definitions. If $\lambda\in\nt^*$, then a subspace $V\subseteq\nt$ is called $\lambda$-\emph{isotropic} if $\lambda([x, y])=0$ for all $x$, $y\in V$. Recall that a~subalgebra $\pt\subseteq\nt$ is called a~\emph{polarization} at $\lambda$ if it is a maximal $\lambda$-isotropic subspace. Polarizations play a key role in the construction of the irreducible representation of $U$ corresponding to a given coadjoint orbit \cite{Kirillov2}. It is known that if $\Theta=U.\lambda$ is the coadjoint $U$-orbit of $\lambda$, then $$\dim\Theta=2\cdot\codim_{\nt}\pt.$$

Let $1\leq k\leq n$. It is very convenient to put
\begin{equation*}
\Ro_k=\{(k,s)\in\Phi^+\mid1\leq s<k\},~\Co_k=\{(r,k)\in\Phi^+\mid {j<k\leq n}\}.
\end{equation*}

\defi{The sets $\Ro_k$, $\Co_k$ are called the $k$th \emph{row} and the
$k$th \emph{column} of $\Phi^+$ respectively. We will write $\row(\alpha)=k$ (resp. $\col(\alpha)=k$) if $\alpha\in\Ro_k$ (resp. $\alpha\in\Co_k$). Note that if $D\subseteq\Phi^+$ is a rook placement, then
\begin{equation}\label{formula:Ro_Co_RP}
|D\cap\Ro_k|\leq1\text{ and }|D\cap\Co_k|\leq1\text{ for all }1\leq k\leq n.
\end{equation}}

\exam{Let $n=6$. On the picture below boxes from $\Ro_5\cup\Co_2$ are grey.
\begin{equation*}
\mymatrix{\pho& \pho& \pho& \pho& \pho& \pho\\
\Top{2pt}\Rt{2pt}\pho& \pho& \pho& \pho& \pho& \pho\\
\pho& \Top{2pt}\Rt{2pt}\gray\pho& \pho& \pho& \pho& \pho\\
\pho& \gray\pho& \Top{2pt}\Rt{2pt}\pho& \pho& \pho& \pho\\
\gray\pho& \gray\gray\pho& \gray\pho& \Top{2pt}\Rt{2pt}\gray\pho& \pho& \pho\\
\pho& \gray\pho& \pho& \pho& \Top{2pt}\Rt{2pt}\pho& \pho\\}
\end{equation*}}

\newpage Let $D$ be a rook placement. Suppose $D=\{\beta_1,\ldots,\beta_s\}$ and denote $$i_r=\row(\beta_r),~j_r=\col(\beta_r),~1\leq r\leq s.$$ Without loss of generality, $j_1<\ldots<j_s$. Put also $j_0=0$, $\Mo_{j_0}=\varnothing$ and
\begin{equation*}
\begin{split}
\Mo_{j_r}&=\left\{(i_r,q)\in\Ro_{i_r}\mid q>j_r\text{ and $(q,j_r)\notin\bigcup\nolimits_{l=0}^{r-1}\Mo_{j_l}$}\right\},\\
\Po_{j_r}&=\{(p,j_r)\in\Co_{j_r}\mid p<i_r\text{ and }(i_r,p)\in\Mo_{j_r}\}\text{ for all }1\leq r\leq s.
\end{split}
\end{equation*}
We put also $\Mo=\bigcup\nolimits_{r=1}^s\Mo_{j_r}$ and $\Po=\bigcup\nolimits_{r=1}^s\Po_{j_r}$. Clearly, $|\Mo_{j_r}|=|\Po_{j_r}|$ for all $r$.

\exam{Let $n=8$ and $D$ be as in Example~\ref{exam:RP1}. Then
\begin{equation*}
\begin{split}
&\Mo_1=\{(3,2)\},~\Mo_2=\{(6,4),~(6,5)\},~\Mo_3=\{(7,4),~(7,5),~(7,6)\},~\Mo_4=\Mo_6=\varnothing,\\
&\Po_1=\{(2,1)\},~\Po_2=\{(4,2),~(5,2)\},~\Po_3=\{(4,3),~(5,3),~(6,3)\},~\Po_4=\Po_6=\varnothing.
\end{split}
\end{equation*}
On the picture below boxes from $\Mo$ are marked by $-$'s and boxes from $\Po$ are marked by $+$'s.
\begin{equation*}
\mymatrix{ \pho& \pho& \pho& \pho& \pho& \pho& \pho& \pho\\
\Top{2pt}\Rt{2pt}\fho+& \pho& \pho& \pho& \pho& \pho& \pho& \pho\\
\otimes& \Top{2pt}\Rt{2pt}\fho-& \pho& \pho& \pho& \pho& \pho& \pho\\
\pho& \fho+& \Top{2pt}\Rt{2pt}\fho+& \pho& \pho& \pho& \pho& \pho\\
\pho& \fho+& \fho+& \Top{2pt}\Rt{2pt}\otimes& \pho& \pho& \pho& \pho\\
\pho& \otimes& \fho+& \fho-& \Top{2pt}\Rt{2pt}\fho-& \pho& \pho& \pho\\
\pho& \pho& \otimes& \fho-& \fho-& \Top{2pt}\Rt{2pt}\fho-& \pho& \pho\\
\pho& \pho& \pho& \pho& \pho& \otimes& \Top{2pt}\Rt{2pt}\pho& \pho\\
}\end{equation*}}

Denote by $\pt$ the subspace of $\nt$ spanned by all $e_{j,i}$, $(i,j)\in\Phi^+\setminus\Mo$. Recall the definition of $w$ from~(\ref{formula:w_D}): if $D=\{(i_1,j_1),\ldots,(i_s,j_s)\}$, then $$w=(i_1,j_1)\ldots(i_s,j_s)\in S_n.$$ Let $l(w)$ be the length of $w$ in $S_n$. Our second result is as follows (cf. \cite[Theorems 1.1, 1.2]{Panov}, \cite[Theorem 1.1, 1.2]{Ignatyev3}, \cite[Theorem 0.2]{Ignatyev4}, \cite[Proposition 4.1]{Ignatyev1}, \cite[Theorem 3.1]{Ignatyev2}).

\mtheo{Let $D$ be a rook placement and $\xi\colon D\to\Cp^{\times}$ a map. Then \label{mtheo:dim}
\begin{equation*}
\begin{split}
&\mathrm{i})~\pt\subseteq\nt\text{ is a polarization at }f_{D,\xi},\\
&\mathrm{ii})~\dim\Theta_{D,\xi}=2\cdot|\Mo|\leq l(w)-|D|,\\
&\mathrm{iii})~\dim\Omega_D=2\cdot|\Mo|+|D|\leq l(w).
\end{split}
\end{equation*}}
The proof is given below, see Propositions~\ref{prop:pol_subalg}, \ref{prop:pol_max_iso}, \ref{prop:dim_Theta_leq} and Corollaries~\ref{coro:dim_Theta_Mo}, \ref{coro:dim_Omega_leq}.
\nota{Note that if $i>j$, then $w(j)=i$ if and only if $(i,j)\in D$. It follows from \cite{Proctor} that if $D$, $D'$ are rook placements, $w$, $w'$ are the corresponding permutations and $w\leq_B w'$, then $D\leq D'$. But the converse statement is not true in general. It is true if $D$ and $D'$ are orthogonal subsets \cite[Theorem 1.10]{Ignatyev1}. On the other hand, let $n=4$, $D=\{(2,1),(3,2),(4,3)\}$ and $D'=\{(3,1),(4,2)\}$. One can easily see that $D\leq D'$, but $w_D$ and $w_{D'}$ are incomparable in the Bruhat--Chevalley order.}

\sst Let $D$ be a rook placement and $\xi\colon D\to\Cp^{\times}$ a map. Recall the definition of $\pt$ from {Sub\-sec\-tion}~\ref{sst:dim_pol_formulir}. In this Subsection, we prove that $\pt$ is a polarization at $f_{D,\xi}$, cf. \cite[Theorem 1.1]{Panov} and \cite[Theorem 1.1]{Ignatyev3}.

\newpage\propp{The subspace $\pt$ is a subalgebra\label{prop:pol_subalg} of $\nt$.}{By definition, $\pt$ is spanned by $e_{j,i}$, where $(i,j)\in\Phi^+\setminus\Mo$. Suppose that $\pt$ is not a subalgebra. Then there exist $j<k<i$ such that $(i,k)$ and $(k,j)$ belong to $\Phi^+\setminus\Mo$, but $(i,j)$ belongs to $\Mo$. We call such a triple $\{i,k,j\}$ \emph{non-admissible}. We can assume without loss of generality that $j$ is minimal among all such triples.

Since $(i,j)\in\Mo$, there exists $s<j$ such that $(i,s)\in D$ and $(i,j)\in\Mo_s$. But $(i,k)\notin\Mo$, so $(k,s)\in\Mo$. Further, $(i,j)\in\Mo_s$ means that $(j,s)\in\Po_s$, hence $(j,s)\notin\Mo$. Thus, $\{k,j,s\}$ is a non-admissible triple. This contradicts the choice of $j$. The proof is complete.}

\propp{The subspace $\pt$ is a maximal $f_{D,\xi}$-isotropic\label{prop:pol_max_iso} subspace of $\nt$.}{Denote $f=f_{D,\xi}$. If $$x=\sum_{(i,j)\in\Phi^+}x_{i,j}e_{j,i},~y=\sum_{(r,s)\in\Phi^+}y_{r,s}e_{s,r}$$ and $f([x,y])\neq0$, then there exist $j<k<i$ such that $(i,j)\in D$ and, say, $x_{i,k}\neq0$, $y_{k,j}\neq0$. If $x,y\in\pt$, then $(i,k),(k,j)\notin\Mo$. This contradicts the definition of $\Mo$. We conclude that $f([x,y])=0$ for all $x,y\in\pt$, hence $\pt$ is an isotropic subspace of $\nt$.

Now, assume that $$x=\sum_{(i,j)\in\Phi^+}x_{i,j}e_{j,i}$$ does not belong to $\pt$, but $\pt+\Cp x$ is an isotropic subspace of $\nt$. Since $x\notin\pt$, there exists a root $(i,j)\in\Mo$ such that $x_{i,j}\neq0$. We can assume that $j$ is minimal among all such roots. It follows from $(i,j)\in\Mo$ that there exists $(i,s)\in D$ such that $(i,j)\in\Mo_s$. Note that $(j,s)\in\Po_s$, hence $e_{s,j}\in\pt$.

We claim that $f([e_{s,j},x])\neq0$. Indeed,
\begin{equation*}
\begin{split}
[e_{s,j},x]&=x_{i,j}[e_{s,j},e_{j,i}]+\sum_{(a,b)\in\Phi^+,~b\geq j,~a\neq i}x_{a,b}[e_{s,j},e_{b,a}]\\
&=x_{i,j}e_{s,i}+\sum_{(a,j)\in\Phi^+,~a\neq i}x_{a,j}[e_{s,j},e_{j,a}]\\
&=x_{i,j}e_{s,i}+\sum_{(a,j)\in\Phi^+,~a\neq i}x_{a,j}e_{s,a}.
\end{split}
\end{equation*}
If $a\neq i$, then $(a,s)\notin D$, so $f(e_{s,a})=0$. Thus, $f([e_{s,j},x])=\xi(i,j)x_{i,j}\neq0$. This contradiction shows that $\pt$ is a maximal isotropic subspace of $\nt$, as required.}

\corop{The dimension of the orbit $\Theta_{D,\xi}$ does not depend on $\xi$ and\label{coro:dim_Theta_Mo} equals $$\dim\Theta_{D,\xi}=2\cdot|\Mo|.$$}{Since $\pt$ is a maximal $f_{D,\xi}$-isotropic subspace, $\dim\Theta_{D,\xi}$ equals $2\cdot\codim\pt.$}

\sst Recall the definition of $w\in S_n$ from {Sub\-sec\-tion}~\ref{sst:dim_pol_formulir}. In this Subsection, we show that $\dim\Theta_{D,\xi}\leq l(w)-|D|$ and, consequently, $\dim\Omega_D\leq l(w)$, cf. \cite[Theorem 1.2]{Panov}, \cite[Theorem 1.2]{Ignatyev3}, \cite[Theorem 0.2]{Ignatyev4}, \cite[Proposition 4.1]{Ignatyev1} and \cite[Theorem 3.1]{Ignatyev2}.

\propp{Let $D$ be a rook placement and $\xi\colon D\to\Cp$ be a map. Then\label{prop:dim_Theta_leq}
\begin{equation*}
\dim\Theta_{D,\xi}\leq l(w)-|D|.
\end{equation*}}{Put $\Theta=\Theta_{D,\xi}$. We will proceed by induction on $n$ (the base is evident). If $D\cap\Co_1=\varnothing$, then the required inequality holds by the inductive hypothesis, so we may assume that $D\cap\Co_1=\{(i_1,1)\}\neq\varnothing$. Suppose that $(i_1,1),(i_2,i_1),\ldots,(i_k,i_{k-1})\in D$ and $D\cap\Ro_{i_k}=\varnothing$. Denote $i_0=1$ and consider the set $\wt\Phi=\pm\wt\Phi^+$, where $$\wt\Phi^+=\Phi^+\setminus\left(\bigcup_{j=0}^k(\Ro_{i_j}\cup\Co_{i_j})\right).$$

We can identify $\wt\Phi$ with the root system $A_{n-k-1}$ by the obvious rule. Put also $\wt D=D\cap\wt\Phi^+$ and $$\wt w=\prod_{(i,j)\in\wt D}(i,j).$$ Let $\wt l(\wt w)$ be the length of $\wt w$ in the Weyl group of $\wt\Phi$. Denote by $\wt U$ the subgroup of $U$ generated by $x_{j,i}(\alpha_{j,i})$, $(i,j)\in\wt\Phi^+$, $\alpha_{j,i}\in\Cp$. Denote also by $\wt\nt$ the subalgebra of $\nt$ generated by $e_{j,i}$, $(i,j)\in\wt\Phi^+$. Finally, denote by $\wt\Theta$ the $\wt U$-orbit of the linear form $$\wt f=\sum_{(i,j)\in\wt D}\xi(i,j)e_{i,j}\in\wt\nt^*.$$

By the inductive assumption, $\dim\wt\Theta\leq\wt l(\wt w)-|\wt D|$. Let $\wt\pt$ be the polarization of $\wt\nt$ at $\wt f$ constructed by the rule described in Subsection~\ref{sst:dim_pol_formulir}. Clearly, $$\pt=\wt\pt\oplus\bigoplus_{(i,j)\in\left(\Phi^+\setminus\wt\Phi^+\right)\setminus\Mo'}e_{j,i},$$ where $$\Mo'=\bigcup_{j=0}^{k-1}\Mo_{i_j}=\bigcup_{r=1}^k\left(\bigcup_{j=i_{r-1}+1}^{i_r-1}\{(i_r,j)\}\right).$$ Using Corollary~\ref{coro:dim_Theta_Mo}, we obtain
\begin{equation}
\begin{split}
\dim\Theta&=2\cdot|\Mo|=\dim\wt\Theta+2\cdot|\Mo'|\\
&=\dim\wt\Theta+2\cdot\sum_{r=1}^k(i_r-i_{r-1}-1)=\dim\wt\Theta+2i_k-2k-2.
\end{split}\label{formula:dim_Theta_wt_Theta}
\end{equation}

Our goal now is to compare $l(w)$ with $\wt l(\wt w)$. Recall that a pair $\{a,b\}$ is called an \emph{inversion} in $w$ if $a>b$ and $w^{-1}(a)<w^{-1}(b)$. It is well-known that the length of $w$ equals the number of inversions in~$w$. Denote $X=\{i_0=1,i_1,\ldots,i_k\}$, then $$l(w)=\wt l(\wt w)+\#\{\{a,b\}\mid\{a,b\}\text{ is an inversion in }w\text{ and }\{a,b\}\cap X\neq\varnothing\}.$$
Note that $w$ has the form $$w=\begin{pmatrix}
1&\ldots&i_1&\ldots&i_2&\ldots&i_{k-2}&\ldots&i_{k-1}&\ldots&i_k&\ldots\\
i_1&\ldots&i_2&\ldots&i_3&\ldots&i_{k-1}&\ldots&i_k&\ldots&1&\ldots\\
\end{pmatrix}.$$
In particular, if $w^{-1}(a)<i_k$, then $\{a,1\}$ is an inversion in $w$, hence $$l(w)\geq\wt l(\wt w)+i_k-1.$$

Furthermore, put $Y_r=\{i_{r-1}+1,i_{r-1}+2,\ldots,i_r-1\}$, $1\leq r\leq k$. Pick a number $y\in Y_k$. If $w^{-1}(y)>i_{k-1}$, then $\{i_k,y\}$ is an inversion in $w$. On the other hand, if $w^{-1}(y)<i_{k-2}$, then $\{y,i_{k-1}\}$ is an inversion in $w$. Finally, suppose $i_{k-2}<w^{-1}(y)<i_{k-1}$. Then there exist $y_1<y_2<\ldots<y_s$ such that $w(y_1)=y_2$, $w(y_2)=y_3$, $\ldots$, $w(y_{s-1})=y_s$, $w(y_s)=y_1$ and $y=y_i$ for some $i$, see the beginning of the proof of Lemma~\ref{lemm:Omega_union_Theta}. In this case, $\{i_k,y_1\}$ is an inversion in $w$, because $i_k>y>y_1$. Thus, $$l(w)\geq\wt l(\wt w)+i_k-1+|Y_k|.$$

Arguing by the similar way, we obtain
\begin{equation*}
\begin{split}
l(w)&\geq\wt l(\wt w)+i_k-1+\sum_{r=1}^k|Y_r|=\wt l(\wt w)+i_k-1+\sum_{r=1}^k(i_r-i_{r-1}-1)\\
&=\wt l(\wt w)+i_k-1+(i_k-1)-k=\wt l(\wt w)+2i_k-k-2.
\end{split}
\end{equation*}
Using (\ref{formula:dim_Theta_wt_Theta}) and the inductive assumption, we conclude that
\begin{equation*}
\begin{split}
\dim\Theta&=\dim\wt\Theta+2_i-2k-2\leq\wt l(\wt w)-|\wt D|+2_k-2k-2\\
&=(\wt l(\wt w)+2i_k-k-2)-(|\wt D|+k)\leq l(w)-|D|.
\end{split}
\end{equation*}
The proof is complete.}

\nota{Suppose $D$ is an orthogonal subset of $A_{n-1}^+$ and $\xi\colon D\to\Cp$ is a map. In this case, $\dim\Theta_{D,\xi}=l(w)-|D|$, see \cite[Theorem~1.2]{Panov}. In general, this is not true. For example, let $n=6$ and $D=\{(3,1),(5,2),(4,3),(6,4)\}$:
\begin{equation*}
\mymatrix{\pho& \pho& \pho& \pho& \pho& \pho\\
\Top{2pt}\Rt{2pt}+& \pho& \pho& \pho& \pho& \pho\\
\otimes& \Top{2pt}\Rt{2pt}-& \pho& \pho& \pho& \pho\\
\pho& +& \Top{2pt}\Rt{2pt}\otimes& \pho& \pho& \pho\\
\pho& \otimes& \pho& \Top{2pt}\Rt{2pt}-& \pho& \pho\\
\pho& \pho& \pho& \otimes& \Top{2pt}\Rt{2pt}\pho& \pho\\}
\end{equation*}
Then, by Corollary~\ref{coro:dim_Theta_Mo}, $\dim\Omega_{D,\xi}=2\cdot|\Mo|=4$. At the contrary, $w=\begin{pmatrix}1&2&3&4&5&6\\3&5&4&6&2&1\end{pmatrix}$, hence $l(w)=10$ and $l(w)-|D|=10-4=6>4$.}

\corop{Let $D$ be a rook placement. Then\label{coro:dim_Omega_leq} $$\dim\Omega_D=2\cdot|\Mo|+|D|\leq l(w)$$.}{Denote by $\xi_0$ the map from $D$ to $\Cp^{\times}$ such that $\xi_0(i,j)=1$ for all $(i,j)\in D$, and $\Theta_0=\Theta_{D,\xi_0}$, so $f_D\in\Theta_0$. Let $Z$ be the stabilizer of $f_D$ in $B$ (cf. the proof of \cite[Proposition 4.1]{Ignatyev1}). Then $$\dim\Omega_D=\dim B-\dim Z.$$ Recall $B=U\rtimes T$. Suppose $g=ut\in Z$, where $u\in U$, $t\in T$, then $g.f_D=u.(t.f_D)=u.f_{D,\xi}$, where $\xi(i,j)=t_{i,i}/t_{j,j}$ (cf. the proof of Lemma~\ref{lemm:Omega_union_Theta}).

Hence $g.f_D\in\Theta_{D,\xi}$. But \cite[Theorem 1]{Andre1} shows that $\Theta_{D,\xi_1}=\Theta_{D,\xi_2}$ if and only if $\xi_1=\xi_2$, so $\xi=\xi_0$ and, consequently, $g.f_D=f_D$. This means that map $$Z\to Z_U\times Z_T\colon ut\mapsto(u,t)$$ is an isomorphism of algebraic varieties, where $Z_U$ (resp. $Z_T$) is the stabilizer of $f_D$ in $U$ (resp. in $T$). Thus, $\dim Z=\dim Z_U+\dim Z_T$.

But $t\in T$ belongs to $Z_T$ if and only if $t_{i,i}/t_{j,j}=1$ for all $(i,j)\in D$, hence $$\dim Z_T=\dim T-|D|=n-|D|.$$
On the other hand, $\dim\Theta_0=2\cdot|\Mo|\leq l(w)-|D|$ by Proposition~\ref{prop:dim_Theta_leq} and Corollary~\ref{coro:dim_Theta_Mo}, thus
\begin{equation*}
\begin{split}
\dim Z_U&=\dim U-\dim\Theta_0\\
&=\dim B-n-2\cdot|\Mo|\geq\dim B-n-l(w)+|D|.
\end{split}
\end{equation*}
We see that $\dim Z=\dim B-2\cdot|\Mo|-|D|\geq\dim B-l(w)$, and so $$\dim\Omega_D=2\cdot|\Mo|+|D|\leq l(w),$$ as required. The proof is complete.}

\sect{The covering relation of the set of rook placements}\label{sect:cov_rel}

\sst\label{sst:cov_rel_formulir} Our third result describes the covering relation of the set of rook placements with the partial order $\leq$ defined above. Denote this set by $\Ro$. To a given $D\in\Ro$, we will describe the set $$L(D)=\{D'\in\Ro\mid D'<D\text{ and there are no }D''\in\Ro\text{ such that }D'<D''<D\}$$ in the spirit of \cite{Ignatyev1}.

There exists a natural partial order on $\Phi^+$: $\alpha\leq\beta$ if $\beta-\alpha$ is a sum of positive roots. In other words, $(a,b)\leq(c,d)$ if $a\leq c$ and $b\geq d$. Let $D$ be a rook placement. Denote by $\wt M(D)$ the set of minimal elements of $D$. Put
\begin{equation*}
\begin{split}
M(D)&=\{(i,j)\in\wt M(D)\mid D\cap\Ro_k\neq\varnothing\text{ and }D\cap\Co_k\neq\varnothing\text{ for all }j<k<i\},\\
N^-(D)&=\{D_{(i,j)}^-,(i,j)\in M(D)\},
\end{split}
\end{equation*}
where $D_{(i,j)}^-=D\setminus\{(i,j)\}$.

Let $(i,j)\in D$. Denote $$m=\min\{k\mid j<k<i\text{ and }D\cap\Co_k=\varnothing\}.$$
Suppose $m$ exists. Further, suppose there are no $(p,q)\in D$ such that $(i,j)>(p,q)$ and $(i,j)\not>(p,q)$. Then put $$D_{(i,j)}^{\to}=(D\setminus\{(i,j)\})\cup\{(i,m)\}.$$
Similarly, suppose $$m=\max\{k\mid j<k<i\text{ and }D\cap\Ro_k=\varnothing\}$$ exists and there are no $(p,q)\in D$ such that $(i,j)>(p,q)$ and $(m,j)\not>(p,q)$. Then denote $$D_{(i,j)}^{\uparrow}=(D\setminus\{(i,j)\})\cup\{(m,j)\}.$$ Now, denote by $B_{(i,j)}(D)$ the set of $(\alpha,\beta)\in D$ such that $(\alpha,\beta)>(i,j)$ and there are no $(p,q)\in D$ satisfying $(i,j)<(p,q)<(\alpha,\beta)$. If $(\alpha,\beta)\in B_{(i,j)}(D)$, then put $$D_{(i,j)}^{(\alpha,\beta)}=(D\setminus\{(i,j),(\alpha,\beta)\})\cup\{(i,\beta),(\alpha,j)\}.$$
We set $$N^0(D)=\bigcup_{(i,j)\in D}\left\{D_{(i,j)}^{\to},D_{(i,j)}^{\uparrow}\right\}\cup\bigcup_{(i,j)\in D}\left\{D_{(i,j)}^{(\alpha,\beta)},(\alpha,\beta)\in B_{(i,j)}(D)\right\}.$$

\exam{Let $n=8$ and $D=\{(3,1),(6,2),(7,3),(5,4),(8,6)\}$, as in Example~\ref{exam:RP1}. Then $\wt M(D)=M(D)=\{(3,1),(5,4),(8,6)\}$ and $(6,2)\in B_{(5,4)}(D)$. On the picture below we draw $D$, $D_{(5,4)}^{(6,2)}$, $D_{(6,2)}^{\uparrow}$ and $D_{(7,3)}^{\to}$.
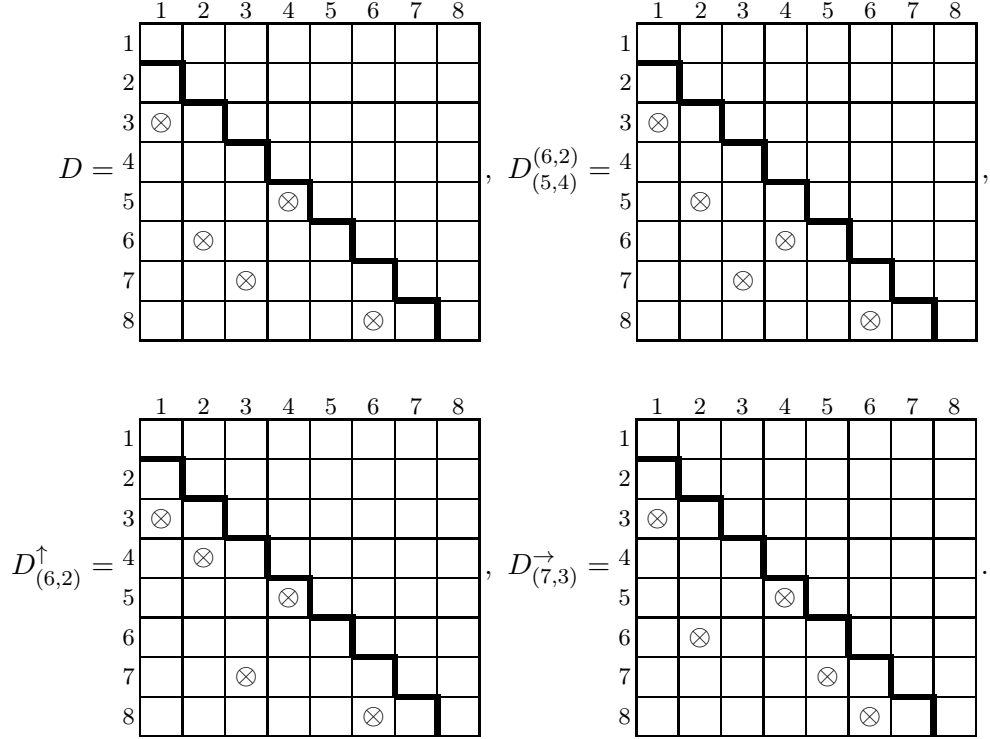
\begin{equation*}
\begin{split}
D&=\mymatrix{ \pho& \pho& \pho& \pho& \pho& \pho& \pho& \pho\\
\Top{2pt}\Rt{2pt}\pho& \pho& \pho& \pho& \pho& \pho& \pho& \pho\\
\otimes& \Top{2pt}\Rt{2pt}\pho& \pho& \pho& \pho& \pho& \pho& \pho\\
\pho& \pho& \Top{2pt}\Rt{2pt}\pho& \pho& \pho& \pho& \pho& \pho\\
\pho& \pho& \pho& \Top{2pt}\Rt{2pt}\otimes& \pho& \pho& \pho& \pho\\
\pho& \otimes& \pho& \pho& \Top{2pt}\Rt{2pt}\pho& \pho& \pho& \pho\\
\pho& \pho& \otimes& \pho& \pho& \Top{2pt}\Rt{2pt}\pho& \pho& \pho\\
\pho& \pho& \pho& \pho& \pho& \otimes& \Top{2pt}\Rt{2pt}\pho& \pho\\}\,,~
D_{(5,4)}^{(6,2)}=\mymatrix{ \pho& \pho& \pho& \pho& \pho& \pho& \pho& \pho\\
\Top{2pt}\Rt{2pt}\pho& \pho& \pho& \pho& \pho& \pho& \pho& \pho\\
\otimes& \Top{2pt}\Rt{2pt}\pho& \pho& \pho& \pho& \pho& \pho& \pho\\
\pho& \pho& \Top{2pt}\Rt{2pt}\pho& \pho& \pho& \pho& \pho& \pho\\
\pho& \otimes& \pho& \Top{2pt}\Rt{2pt}\pho& \pho& \pho& \pho& \pho\\
\pho& \pho& \pho& \otimes& \Top{2pt}\Rt{2pt}\pho& \pho& \pho& \pho\\
\pho& \pho& \otimes& \pho& \pho& \Top{2pt}\Rt{2pt}\pho& \pho& \pho\\
\pho& \pho& \pho& \pho& \pho& \otimes& \Top{2pt}\Rt{2pt}\pho& \pho\\}\,,\\\\
D_{(6,2)}^{\uparrow}&=\mymatrix{ \pho& \pho& \pho& \pho& \pho& \pho& \pho& \pho\\
\Top{2pt}\Rt{2pt}\pho& \pho& \pho& \pho& \pho& \pho& \pho& \pho\\
\otimes& \Top{2pt}\Rt{2pt}\pho& \pho& \pho& \pho& \pho& \pho& \pho\\
\pho& \otimes& \Top{2pt}\Rt{2pt}\pho& \pho& \pho& \pho& \pho& \pho\\
\pho& \pho& \pho& \Top{2pt}\Rt{2pt}\otimes& \pho& \pho& \pho& \pho\\
\pho& \pho& \pho& \pho& \Top{2pt}\Rt{2pt}\pho& \pho& \pho& \pho\\
\pho& \pho& \otimes& \pho& \pho& \Top{2pt}\Rt{2pt}\pho& \pho& \pho\\
\pho& \pho& \pho& \pho& \pho& \otimes& \Top{2pt}\Rt{2pt}\pho& \pho\\}\,,~
D_{(7,3)}^{\to}=\mymatrix{ \pho& \pho& \pho& \pho& \pho& \pho& \pho& \pho\\
\Top{2pt}\Rt{2pt}\pho& \pho& \pho& \pho& \pho& \pho& \pho& \pho\\
\otimes& \Top{2pt}\Rt{2pt}\pho& \pho& \pho& \pho& \pho& \pho& \pho\\
\pho& \pho& \Top{2pt}\Rt{2pt}\pho& \pho& \pho& \pho& \pho& \pho\\
\pho& \pho& \pho& \Top{2pt}\Rt{2pt}\otimes& \pho& \pho& \pho& \pho\\
\pho& \otimes& \pho& \pho& \Top{2pt}\Rt{2pt}\pho& \pho& \pho& \pho\\
\pho& \pho& \pho& \pho& \otimes& \Top{2pt}\Rt{2pt}\pho& \pho& \pho\\
\pho& \pho& \pho& \pho& \pho& \otimes& \Top{2pt}\Rt{2pt}\pho& \pho\\}\,.
\end{split}
\end{equation*}}

Now, suppose $(i,j)\in D$ and $(\alpha,\beta)\in\Zp\times\Zp$. Assume that $i>\beta\geq\alpha>j$, $D\cap\Ro_{\alpha}=D\cap\Co_{\beta}=\varnothing$, $D\cap\Ro_k\neq\varnothing$, $D\cap\Co_k\neq\varnothing$ for all $\alpha<k<\beta$ and if $(p,q)\in D$, $(i,j)>(p,q)$, $(\alpha,j)\not>(p,q)$, then $(i,\beta)>(p,q)$. Further, assume that if $\alpha\neq\beta$, then $D\cap\Ro_{\beta}\neq\varnothing$ and $D\cap\Co_{\alpha}\neq\varnothing$. Denote the set of such pairs $(\alpha,\beta)$ by $C_{(i,j)}(D)$. If $(\alpha,\beta)\in C_{(i,j)}(D)$, then we put $$D_{(i,j)}^{\alpha,\beta}=(D\setminus\{(i,j)\})\cup\{(i,\beta),(\alpha,j)\}.$$
We also set $$N^+(D)=\bigcup_{(i,j)\in D}\left\{D_{(i,j)}^{\alpha,\beta},(\alpha,\beta)\in C_{(i,j)}(D)\right\}.$$

\exam{Let $n=6$ and $D=\{(4,1),(6,2),(5,4)\}$. Then $(3,3)\in C_{(6,2)}(D)$. On the picture below we draw $D$ and $D_{(6,2)}^{3,3}$.
\begin{equation*}
D=\mymatrix{ \pho& \pho& \pho& \pho& \pho& \pho\\
\Top{2pt}\Rt{2pt}\pho& \pho& \pho& \pho& \pho& \pho\\
\pho& \Top{2pt}\Rt{2pt}\pho& \pho& \pho& \pho& \pho\\
\otimes& \pho& \Top{2pt}\Rt{2pt}\pho& \pho& \pho& \pho\\
\pho& \pho& \pho& \Top{2pt}\Rt{2pt}\otimes& \pho& \pho\\
\pho& \otimes& \pho& \pho& \Top{2pt}\Rt{2pt}\pho& \pho\\}\,,~
D_{(6,2)}^{3,3}=\mymatrix{ \pho& \pho& \pho& \pho& \pho& \pho\\
\Top{2pt}\Rt{2pt}\pho& \pho& \pho& \pho& \pho& \pho\\
\pho& \Top{2pt}\Rt{2pt}\otimes& \pho& \pho& \pho& \pho\\
\otimes& \pho& \Top{2pt}\Rt{2pt}\pho& \pho& \pho& \pho\\
\pho& \pho& \pho& \Top{2pt}\Rt{2pt}\otimes& \pho& \pho\\
\pho& \pho& \otimes& \pho& \Top{2pt}\Rt{2pt}\pho& \pho\\}\,.
\end{equation*}}

Finally, we denote
\begin{equation*}
N(D)=N^-(D)\cup N^0(D)\cup N^+(D)
\end{equation*}
(cf. \cite[Subsection~2.2]{Ignatyev1} and \cite[Subsections~3.7--3.14]{Melnikov2}).

\newpage\mtheo{Let $D\in\Ro$ be a rook placement.\label{mtheo:cov_rel} Then $L(D)=N(D)$.}

The proof is given below. In~Subsection~\ref{sst:cov_ness}, we check that $N(D)\subseteq L(D)$. In Subsection~\ref{sst:cov_suff}, we prove that $L(D)\subseteq N(D)$, so $N(D)=L(D)$.

\sst\label{sst:cov_ness} Let $D\in\Ro$ be a rook placement. Recall the definitions of $N^-(D)$, $N^0(D)$, $N^+(D)$,~$N(S)$ and~$L(D)$ from Subsection~\ref{sst:cov_rel_formulir}. Note that
\begin{equation}
\begin{split}
&|T|=|D|\pm1\text{ for all }T\in N^{\pm}(D),\\
&|T|=|D|\text{ for all }T\in N^0(D).
\end{split}\label{formula:dlina_N_pm_0}
\end{equation}
Theorem~\ref{mtheo:cov_rel} claims that $L(D)=N(D)$. To prove this, we need some more notation. Namely, put
\begin{equation*}
\begin{split}
L^-(D)&=\{T\in\Ro\mid T<D,~|T|<|D|\text{ and if }T\leq S<D\text{ then }S=T\},\\
L^0(D)&=\{T\in\Ro\mid T<D,~|T|=|D|\text{ and if }T\leq S<D\text{ then }S=T\},\\
L^+(D)&=\{T\in\Ro\mid T<D,~|T|>|D|\text{ and if }T\leq S<D\text{ then }S=T\},\\
\wt L^-(D)&=\{T\in\Ro\mid T<D,~|T|<|D|\text{ and if }T\leq S<D,~|S|<|D|,\text{ then }S=T\},\\
\wt L^0(D)&=\{T\in\Ro\mid T<D,~|T|=|D|\text{ and if }T\leq S<D,~|S|\leq|D|,\text{ then }S=T\},\\
\wt N^-(D)&=\{D_{(i,j)}^-,(i,j)\in\wt M(D)\}.\\
\end{split}
\end{equation*}
Clearly, $N^-(D)\subseteq\wt N^-(D)$, $L^-(D)\subseteq\wt L^-(D)$, $L^0(D)\subseteq\wt L^0(D)$ and $L(D)=L^-(D)\cup L^0(D)\cup L^+(D)$. Our strategy is to prove that $$N^-(D)=L^-(D),~N^0(D)=L^0(D),~N^+(D)=L^+(D).$$
In this Subsection, we will prove that $$\wt N^-(D)\subseteq \wt L^-(D),~N^0(D)\subseteq L^0(D),~N^+(D)\subseteq L^+(D).$$

\defi{A subset $Y\subseteq\Phi^+$ is called an \emph{ideal} if it follows from $(a,b)\in Y$, $(c,d)\in\Phi^+$, $(c,d)>(a,b)$ that $(c,d)\in Y$.
One can easily see that if $D,D'\in\Ro$ are rook placements and $Y$ is an ideal, then
\begin{equation}
Y\cap D=Y\cap D'\Longleftrightarrow(R_D)_{i,j}=(R_{D'})_{i,j}\text{ for all }(i,j)\in Y.\label{formula:Y_cap_D}
\end{equation}}

Now we are ready to prove that $\wt N^-(D)\subseteq\wt L^-(D)$ (cf. \cite[Lemma 3.8]{Melnikov2}, \cite[Lemma 3.2]{Ignatyev1}).
\lemmp{Let $D\in\Ro$.\label{lemm:wt_N_minus_in_wt_L_minus} One has $\wt N^-(D)\subseteq\wt L^-(D)$.}{Suppose $T=D_{(i,j)}^-\in\wt N^-(D)$ for some root $(i,j)\in\wt M(D)$. By~(\ref{formula:dlina_N_pm_0}), $|T|=|D|-1<|D|$. Put
$$Y=\{(p,q)\in\Phi^+\mid(p,q)\nleq(i,j)\},$$ then $T=Y\cap T=Y\cap D$. Assume that there exists $S\in\Ro$ such that $T\leq S<D$ and $|S|<|D|$. By~(\ref{formula:Y_cap_D}), $$Y\cap T= Y\cap S= Y\cap D,$$ so $|S|\leq|T|=|D|-1$. Thus, $|S|=|T|$ and $S=T$. This means that $T\in\wt L^-(D)$.}

Second, let us show that $N^0(D)\subseteq L^0(D)$ (cf. \cite[Lemmas 3.11--3.14]{Melnikov2}, \cite[Lemma 3.3]{Ignatyev1}).

\lemmp{Let $D\in\Ro$. One\label{lemm:N_0_in_L_0} has $N^0(D)\subseteq L^0(D)$.}{Let $T\in N^0(D)$. By (\ref{formula:dlina_N_pm_0}), $|T|=|D|$. First, assume $T=D_{(i,j)}^{\uparrow}$ for some $(i,j)\in D$. (The case $T=D_{(i,j)}^{\to}$ is completely similar.) Suppose $T\setminus D=\{(m,j)\}$. Put $Y=Y_0\cup Y_1$
and $\wt Y=\Phi^+\setminus Y$, where
\begin{equation*}
\begin{split}
&Y_0=\{(p,q)\in\Phi^+\mid(p,q)\nleq(i,j)\},\\
&Y_1=\{(p,q)\in\Phi^+\mid(p,q)\leq(m,j)\}.
\end{split}
\end{equation*}
Then $(R_D)_{r,s}=(R_T)_{r,s}$ for all $(r,s)\in
Y$. Note that $Y_0$ and $\wt Y\cup Y_0$ are ideals. For example, let $n=8$, $i=7$, $j=2$, $m=5$. On the picture below
boxes from $Y_0$ are filled by~$0$'s, boxes from $Y_1$ are filled by~$1$'s, and boxes from $\wt Y$ are grey.
\begin{equation*}
\mymatrix{\pho& \pho& \pho& \pho& \pho& \pho& \pho& \pho\\
\Top{2pt}\Rt{2pt}0& \pho& \pho& \pho& \pho& \pho& \pho& \pho\\
\Rt{2pt}0& \Top{2pt}\Rt{2pt}1& \pho& \pho& \pho& \pho& \pho& \pho\\
\Rt{2pt}0& 1& \Top{2pt}\Rt{2pt}1& \pho& \pho& \pho& \pho& \pho\\
\Rt{2pt}0& 1& 1& \Top{2pt}\Rt{2pt}1& \pho& \pho& \pho& \pho\\
\Rt{2pt}0& \Top{2pt}\gray\pho& \Top{2pt}\gray\pho& \Top{2pt}\gray\pho& \Top{2pt}\Rt{2pt}\gray\pho& \pho& \pho& \pho\\
\Rt{2pt}0& \gray\pho& \gray\pho& \gray\pho& \gray\pho& \Top{2pt}\Rt{2pt}\gray\pho& \pho& \pho\\
0& \Top{2pt}0& \Top{2pt}0& \Top{2pt}0& \Top{2pt}0& \Top{2pt}0& \Top{2pt}\Rt{2pt}0& \pho\\
}\end{equation*}

Now, assume that there exists $S\in\Ro$ such that $T\leq S<D$. By (\ref{formula:Y_cap_D}), it is enough to check that $(R_S)_{r,s}=(R_T)_{r,s}$ for all $(r,s)\in\wt Y$, because $(R_D)_{r,s}=(R_T)_{r,s}=(R_S)_{r,s}$ for all $(r,s)\in Y$. Note that $$(R_T)_{r,s}=(R_D)_{r,s}-1$$ for all $(r,s)\in\wt Y$. Moreover, by definition of $D_{i,j}^{\uparrow}$, $D\cap\wt Y=\{(i,j)\}$ and $T\cap\wt Y=\varnothing$.

It follows from $T\leq S<D$ that there exists $(k,j)\in S$ such that $m\leq k\leq i$. We claim that $S\cap(\wt Y\setminus\Ro_i)=\varnothing$. Indeed, assume there exists $(p,q)\in S\cap(\wt Y\setminus \Ro_i)$. Then $m<p<i$. By definition of~$D_{i,j}^{\uparrow}$, $D\cap Y_0\cap\Ro_p\neq\varnothing$. Since $D\cap Y_0=T\cap Y_0=S\cap Y_0$, we obtain $S\cap Y_0\cap\Ro_p\neq\varnothing$. But $(p,q)\in\wt Y$, hence $|S\cap\Ro_p|\geq2$. This contradicts (\ref{formula:Ro_Co_RP}). Thus, $S\cap(\wt Y\setminus\Ro_i)=\varnothing$. In particular, either $k=i$ or $k=m$. If $k=i$, then $$S\cap(\wt Y\cup Y_0)=D\cap(\wt Y\cup Y_0).$$ By (\ref{formula:Y_cap_D}), $S=D$, a contradiction. Hence $k=m$, so $S=T$, as required.

Second, assume that $T=D_{(i,j)}^{(\alpha,\beta)}$ for some $(i,j)\in D$, $(\alpha,\beta)\in B_{(i,j)}(D)$. Suppose $S\in\Ro$ and $T\leq S<D$. Arguing as in step (ii) of the proof of \cite[Lemma 3.3]{Ignatyev1}, we deduce that $S=T$. The proof is complete.}

Finally, we will check that $N^+(D)\subseteq L^+(D)$ (cf. \cite[Lemma 3.4]{Ignatyev1}).

\lemmp{Let $D\in\Ro$. One\label{lemm:N_+_in_L_+} has $N^+(D)\subseteq L^+(D)$.}{Let $T=D_{(i,j)}^{\alpha,\beta}\in N^+(D)$ for some $(i,j)\in D$, $(\alpha,\beta)\in C_{(i,j)}(D)$. By (\ref{formula:dlina_N_pm_0}), $|T|=|D|+1$. Here we put
$Y=Y_0\cup Y_1$ and $\wt Y=\Phi^+\setminus Y$, where
\begin{equation*}
\begin{split}
&Y_0=\{(p,q)\in\Phi^+\mid(p,q)\nleq(i,j)\},\\
&Y_1=\{(p,q)\in\Phi^+\mid(p,q)\leq(\alpha,j)\text{ or
}(p,q)\leq(i,\beta)\}.
\end{split}
\end{equation*}
Then $(R_D)_{r,s}=(R_T)_{r,s}$ for all $(r,s)\in Y$. Note that $Y_0$ and $\wt Y\cup Y_0$ are ideals. For example, let $n=8$, $i=7$, $j=2$, $\alpha=4$, $\beta=5$. On the picture below boxes from $Y_0$ are filled by~$0$'s, boxes
from $Y_1$ are filled by $1$'s, and boxes from $\wt Y$ are grey.
\begin{equation*}
\mymatrix{\pho& \pho& \pho& \pho& \pho& \pho& \pho& \pho\\
\Top{2pt}\Rt{2pt}0& \pho& \pho& \pho& \pho& \pho& \pho& \pho\\
\Rt{2pt}0& \Top{2pt}\Rt{2pt}1& \pho& \pho& \pho& \pho& \pho& \pho\\
\Rt{2pt}0& 1& \Top{2pt}\Rt{2pt}1& \pho& \pho& \pho& \pho& \pho\\
\Rt{2pt}0& \gray\Top{2pt}\pho& \gray\Top{2pt}\pho& \gray\Top{2pt}\Rt{2pt}\pho& \pho& \pho& \pho& \pho\\
\Rt{2pt}0& \gray\pho& \gray\pho& \Rt{2pt}\gray\pho& \Top{2pt}\Rt{2pt}1& \pho& \pho& \pho\\
\Rt{2pt}0& \gray\pho& \gray\pho& \Rt{2pt}\gray\pho& 1& \Top{2pt}\Rt{2pt}1& \pho& \pho\\
0& \Top{2pt}0& \Top{2pt}0& \Top{2pt}0& \Top{2pt}0& \Top{2pt}0& \Top{2pt}\Rt{2pt}0& \pho\\
}\end{equation*}

Now suppose there exists $S\in\Ro$ such that $T\leq S<D$. By (\ref{formula:Y_cap_D}), it is enough to show that $(R_S)_{r,s}=(R_T)_{r,s}$ for all $(r,s)\in\wt Y$. Note that $$(R_T)_{r,s}=(R_D)_{r,s}-1$$ for all $(r,s)\in\wt Y$. Moreover, by definition of $D_{(i,j)}^{\alpha,\beta}$, $D\cap\wt Y=\{(i,j)\}$ and $T\cap\wt Y=\varnothing$.

Since $T\leq S<D$, there exists $(k,j)\in S$ such that $\alpha\leq k\leq i$. If $k=i$, then $$S\cap(\wt Y\cup Y_0)=D\cap(\wt Y\cup Y_0).$$ By (\ref{formula:Y_cap_D}), $S=D$, a contradiction. Hence $\alpha\leq k<i$. Similarly, there exists $(i,l)\in S$ such that $j<l\leq\beta$. We claim that either $l\leq\alpha$ or $l=\beta$. Indeed, assume $\alpha<l<\beta$. By definition of $D_{(i,j)}^{\alpha,\beta}$, $D\cap Y_0\cap\Co_l\neq\varnothing$. Since $D\cap Y_0=T\cap Y_0=S\cap Y_0$, we obtain $S\cap Y_0\cap\Co_l\neq\varnothing$. But $(i,l)\in\wt Y$, hence $|S\cap\Co_l|\geq2$. This contradicts (\ref{formula:Ro_Co_RP}). Thus, either $l\leq\alpha$ or $l=\beta$.

Assume $l=\alpha$. If $k>\alpha$, then $(R_S)_{k,\alpha}>(R_D)_{k,\alpha}$, which contradicts $S<D$, hence $k=\alpha$. If $\alpha<\beta$, then, by definition of $D_{(i,j)}^{\alpha,\beta}$, $$S\cap Y_0\cap\Co_{\alpha}=D\cap Y_0\cap\Co_{\alpha}\neq\varnothing,$$ so $|S\cap\Co_{\alpha}|\geq2$, a contradiction. Hence $\alpha=\beta=k=l$, so $S=T$, as required.

If $l<\alpha$, then $(R_S)_{k,l}>(R_D)_{k,l}$, because $k\geq\alpha$. This contradicts $S<D$, so $l=\beta$. Similarly, $k=\alpha$, hence $$S\cap(\wt Y\cup Y_0)=T\cap(\wt Y\cup Y_0).$$ By (\ref{formula:Y_cap_D}), $S=T$, as required.}

\sst\label{sst:cov_suff} In this Subsection, we will show that $$\wt L^-(D)\subseteq\wt N^-(D),~L^0(D)\subseteq N^0(D),~L^+(D)\subseteq N^+(D).$$
Combining this with the results of the previous Subsection, we see that $$\wt N^-(D)=\wt L^-(D),~N^0(D)=L^0(D),~N^+(D)=L^+(D).$$ Using Lemma~\ref{lemm:wt_N_minus_in_wt_L_minus} and the fact that $\wt N^-(D)=\wt L^-(D)$, one can mimic the proof of \cite[Lemma 3.8]{Ignatyev1} to show that $N^-(D)=L^-(D)$, hence $N(D)=L(D)$. This concludes the proof of Theorem~\ref{mtheo:cov_rel}.

The proofs are much more complicated than the proofs in the previous Subsection. Note that $$D_0=\{(n,1),(n-1,2),\ldots,(n-n_0+1,n_0)\},$$ where $n_0=[n/2]$, is the maximal element of $\Ro$ with respect to the partial order $\leq$ on $\Ro$.

First, we will prove that $\wt L^-(D)=\wt N^-(D)$ (cf. \cite[Lemma 3.5]{Ignatyev1}).

\lemmp{Let $D\in\Ro$. One has $\wt L^-(D)=\wt N^-(D)$.}{By Lemma~\ref{lemm:wt_N_minus_in_wt_L_minus}, it is enough to check that $\wt L^-(D)\subseteq\wt N^-(D)$. We must show that
\begin{equation}
\begin{split}
&\text{if }
T\leq D\text{ and }|T|<|D|,\\
&\text{then there exists }S\in\wt N^-(D)
\text{ such that }T\leq S<D.
\end{split}\label{formula:proof_minus}
\end{equation}
We will proceed by induction on $n$ (for $n=1$, there is nothing to prove). The proof is rather long, so we split it into five
steps.

i) One can easily check that if $D=D_0$, then (\ref{formula:proof_minus}) holds. Therefore, we may also use the second (downward) induction on the partial order $\leq$ on $\Ro$.

ii) Let $D<D_0$, $D=\{(i_1,j_1),\ldots,(i_s,j_s)\}$, as above, $T=\{(p_1,q_1),\ldots,(p_t,q_t)\}$, $p_l>q_l$, $q_l<q_{l+1}$, $T<D$ and $t<s$. Consider the following conditions.
\begin{equation}
\begin{split}
&\text{a) There exists $k\leq n$ such that $i_k=j_k+1$ or $k=|D|$.}\\
&\text{b) There exists $d\leq k$ such that either $q_{d+1}>j_k$ or $d=|T|$.}\\
&\text{c) $i_l>i_{l+1}$ for any $1\leq l\leq k-1$.}\\
&\text{d) $p_l>p_{l+1}$ for any $1\leq l\leq d-1$.}\\
&\text{e) $i_l\geq p_l$ for any $1\leq l\leq d$.}\\
&\text{f) $j_l=l$ for any $1\leq l\leq k$ and $q_l=l$ for any $1\leq l\leq d$.}\\
\end{split}
\label{formula:conditions_minus}
\end{equation}
We claim that
\begin{equation}
\begin{split}
&\text{if (\ref{formula:proof_minus}) holds for all $D$, $T$ satisfying (\ref{formula:conditions_minus}),}\\
&\text{then (\ref{formula:proof_minus}) holds for all $D$, $T\in\Ro$.}
\end{split}
\label{formula:if_conds_then_OK_minus}
\end{equation}

To prove this, we need some more notation. Given $A\in\Ro$, we denote $A_r=A\cap\left(\bigcup_{1\leq l\leq r}\Co_l\right)$. Clearly, to prove (\ref{formula:if_conds_then_OK_minus}), it is enough to show that if (\ref{formula:proof_minus}) holds for all $D$, $T$ satisfying (\ref{formula:conditions_minus}), and $D_r$,~$T_r$ do not satisfy~(\ref{formula:conditions_minus}) for some $r$, then (\ref{formula:proof_minus}) holds for $D$, $T$. We will proceed by induction on $r$ (the base $r=1$ is clear). Evidently, we may assume that $T\in\wt L^-(D)$.

iii) Suppose $1\leq r\leq s$ and $D_r,T_r$ satisfy (\ref{formula:conditions_minus}). To perform the induction step, we must prove that either $D_{r+1},T_{r+1}$ satisfy (\ref{formula:conditions_minus}), or (\ref{formula:proof_minus}) holds for $D$, $T$. This is trivially true if $i_k=k+1$ for $D_r$, so we may assume that $\Co_l\cap D\neq\{(l+1,l)\}$ for all $1\leq l\leq r$.

First, consider the case $\Co_{r+1}\cap D=\varnothing$. (And, consequently, $\Co_l\cap D=\varnothing$ for any $k+1\leq l\leq r+1$.) Given $\alpha=(i,j)\in\Phi^+$ and $1\leq a\leq n$, denote
\begin{equation}
\alpha^a=\begin{cases}
(i,j+1),&\text{if }j\leq a,\\
\alpha,&\text{if }j>a.
\end{cases}
\end{equation}
Given $A\in\Ro$, denote by $[A]_a$ (resp. by $[A]^a$) the subset of $\Phi^+$ defined by $[A]_a=\{\alpha^a,~\alpha\in A\}$ (resp. by $[A]^a=\{\alpha\in\Phi^+\mid\alpha^a\in A\}$). Suppose $p_d>d+1$ for $T_r$. Then $[T]_k$ and $[D]_k$ are rook placements with no rooks in the first column, $[T]_k<[D]_k$ and $t=|[T]_k|<s=|[D]_k|$. (Note that if $\Co_k\cap T\neq\varnothing$, then $\Co_{k+1}\cap T=\varnothing$.) By the inductive assumption on $n$, there exists $S_1\in\wt N^-([D]_k)$ such that $[T]_k\leq S_1$. One can easily see that $S=[S_1]^k\in\wt N^-(D)$ and $T\leq S$, hence we are done.

On the other hand, suppose $p_d=d+1$. (And so $\Co_l\cap T=\varnothing$ for any $d+1\leq l\leq r$.) Put $T'=T\setminus\{(d,d+1)\}$. Then $[T']_k$, $[D]_k$ are rook placements, $[T']_k<[D]_k$ and $t-1=|[T']_k|<s=|[D]_k|$. By the inductive assumption on $n$, there exists $S_1=([D]_k)_{(i,j)}^-\in\wt N^-([D]_k)$ such that $[T']_k\leq S_1$. If $(i,j)\neq(i_d,d)$, then $S=[S_1]^k\in\wt N^-(D)$ and $T\leq S$, hence we are done. Finally, if $(i,j)=(i_d,d)$, then $D'>T'$, where $D'=D\setminus\{(i_d,d)\}$. Thus, $[T']_d$, $[D']_d$ are rook placements, $[T']_d<[D']_d$ and $t-1=|[T']_d|<s-1=|[D']_d|$. By the inductive assumption on $n$, there exists $S_2=([D']_d)_{(a,b)}^-\in\wt N^-([D']_d)$ such that $[T']\leq S_2$. Note that $b>d$, thus $S=[S_2]^d\cup\{(i_d,d)\}\in\wt N^-(D)$ and $T\leq S$, as required.

iv) Second, consider the case the case $\Co_{r+1}\cap D\neq\varnothing$. If $\Co_r\cap D=\varnothing$ (and so $\Co_l\cap D=\varnothing$ for any $k+1\leq l\leq r$), then we can argue as on the previous step. Let $\Co_r\cap D\neq\varnothing$, i.e., $k=r$ for $D_r$. Assume $i_{r+1}>i_r$. Then put $$\wt D=\left(D\setminus\{(i_r,r),(i_{r+1},r+1)\}\right)\cup\{(i_{r+1},r),(i_r,r+1)\}.$$
Evidently, $\wt D>D>T$ and $s=|\wt D|>t=|T|$. By the inductive assumption on $\leq$, there exists $\wt S=(\wt D)_{(i,j)}^-\in\wt N^-(\wt D)$ such that $T\leq \wt S$. Note that $j\geq r+1$. If $j>r+1$, then $S=D_{(i,j)}^-\in\wt N^-(D)$ and $T\leq S$, as required.

If $j=r+1$ and $T\cap\Co_r=\varnothing$, then $S=D_{(i_r,r)}^-\geq T$, as required. If $j=r+1$ and $T\cap\Co_r\neq\varnothing$, then put $D'=D\setminus\{(i_r,r)\}$, $T'=T\setminus\{(p_r,r)\}$. We see that $D'>T'$, $[D']_r$ and $[T']_r$ are rook placements, $[D']_r>[T']_r$, and $t-1=|[T']_r|<s-1=|[D']_r|$. By the inductive assumption on $n$, there exists $S_1=([D']_r)_{(a,b)}^-\in\wt N^-([D']_r)$ such that $[T']_r\leq S_1$. Since $b>r$, we conclude that $S=[S_1]^r\cup\{(i_r,r)\}\in\wt N^-(D)$ and $T\leq S$, so we are done.

On the other hand, assume that $i_{r+1}<i_r$. If $T\cap\Co_{r+1}=\varnothing$, then $D_{r+1}$, $T_{r+1}$ satisfy (\ref{formula:conditions_minus}). If ${T\cap\Co_{r+1}=\{(p,r+1)\}\neq\varnothing}$, but $T\cap\Co_r=\varnothing$, then $T<\wt T_1<D$ and $|\wt T_1|=t<s=|D|$, where $\wt T_1=\left(T\setminus\{(p,r+1)\}\right)\cup\{(p,r)\}$. Hence $T\notin\wt L^-(D)$, a contradiction. Now, suppose $T\cap\Co_r\neq\varnothing$ (so $d=r$ for $T_r$) and $T\cap\Co_{r+1}\neq\varnothing$. If $p_{r+1}>p_r$, then denote
\begin{equation*}
\wt T_2=\begin{cases}
\left(T\setminus\{(p_r,r),(p_{r+1},r+1)\}\right)\cup\{(p_{r+1},r),(p_r,r+1)\},&\text{if }p_r>r+1,\\
\left(T\setminus\{(r+1,r),(p_{r+1},r+1)\}\right)\cup\{(p_{r+1},r)\},&\text{if }p_r=r+1.
\end{cases}
\end{equation*}
One can check that $T<\wt T_2<D$ and $|\wt T_2|\leq t<|D|=s$, so $T\notin\wt L^-(D)$. Thus, (\ref{formula:if_conds_then_OK_minus}) is proved.

v) From now on, we may assume that $D$, $T$ satisfy (\ref{formula:conditions_minus}). Note that $(i_k,k)\in M(D)$. If $i_k>k+1$ (and so $k=s$), then $D\cap\Co_{k+1}=\varnothing$, hence we can argue as on step (iii). Suppose $i_k=k+1$. If $T\cap\Co_k\neq\varnothing$, then $d=k$ and $p_k=i_k=k+1$. Applying the induction hypothesis on $n$ to the rook placements $[D']_k$, $[T']_k$, where $D'=D\setminus\{(k+1,k)\}$, $T'=T\setminus\{(k+1,k)\}$, we conclude that there exists $S_1=([D']_k)_{(i,j)}^-\in\wt N^-([D']_k)$ such that $[T']_k<S_1$. Since $j>k$, $S=[S_1]^k\cup\{(k+1,k)\}\in\wt N^-(D)$ and $T<S$. Finally, if $T\cap\Co_k=\varnothing$, then $S=D_{(k+1,k)}^->T$. The proof is complete.}

Second, we will prove that $L^0(D)=N^0(D)$ (cf. \cite[Lemma 3.6]{Ignatyev1}).

\lemmp{Let $D\in\Ro$. One has $L^0(D)=N^0(D)$.}{By Lemma~\ref{lemm:N_0_in_L_0}, it is enough to check that $L^0(D)\subseteq N^0(D)$. Since $L^0(D)\subseteq\wt L^0(D)$, this is equivalent to $\wt L^0(D)\subseteq N^0(D)$. In other words, we must show that
\begin{equation}
\begin{split}
&\text{if }
T\leq D\text{ and }|T|=|D|,\\
&\text{then there exists }S\in N'(D)=\wt N^-(D)\cup N^0(D)
\text{ such that }T\leq S<D.
\end{split}\label{formula:proof_0}
\end{equation}
We will proceed by induction on $n$ (for $n=1$, there is nothing to prove). For convenience, we split the proof into five steps.

i) One can easily check that if $D=D_0$, then (\ref{formula:proof_minus}) holds. Therefore, we may also use the second (downward) induction on the partial order $\leq$ on $\Ro$.

ii) Let $D<D_0$, $D=\{(i_1,j_1),\ldots,(i_t,j_t)\}$, $T=\{(p_1,q_1),\ldots,(p_t,q_t)\}$ and $T<D$. Consider the following conditions.
\begin{equation}
\begin{split}
&\text{a) There exists $k\leq |D|$ such that $i_l>j_l+1$ for any $1\leq l\leq k-1$.}\\
&\text{b) There exists $d\leq k$ such that $i_l\geq p_l$ for any $1\leq l\leq d$.}\\
&\text{c) If $d<k$, then either $i_k=j_k+1$ or $k=|D|$.}\\
&\text{d) If $d<k$, then either $q_{d+1}>j_k$ or $d=|T|$; if $d=k$, then $p_d=d+1$.}\\
&\text{e) $p_l>p_{l+1}$ for any $1\leq l\leq d-1$.}\\
&\text{f) $j_l=l$ for any $1\leq l\leq k$ and $q_l=l$ for any $1\leq l\leq d$.}\\
\end{split}
\label{formula:conditions_0}
\end{equation}

We claim that
\begin{equation}
\begin{split}
&\text{if (\ref{formula:proof_0}) holds for all $D$, $T$ satisfying (\ref{formula:conditions_0}),}\\
&\text{then (\ref{formula:proof_0}) holds for all $D$, $T\in\Ro$.}
\end{split}
\label{formula:if_conds_then_OK_0}
\end{equation}
Clearly, to prove this, it is enough to show that if (\ref{formula:proof_0}) holds for all $D$, $T$ satisfying (\ref{formula:conditions_0}), and $D_r$,~$T_r$ do not satisfy (\ref{formula:conditions_0}) for some $r$, then (\ref{formula:proof_0}) holds for $D$, $T$. We will proceed by induction on $r$ (the base $r=1$ is clear). Evidently, we may assume that $T\in\wt L^0(D)$.

iii) Suppose $1\leq r\leq s$ and $D_r,T_r$ satisfy (\ref{formula:conditions_0}). To perform the induction step, we must prove that either $D_{r+1},T_{r+1}$ satisfy (\ref{formula:conditions_0}), or (\ref{formula:proof_0}) holds for $D$, $T$. This is trivially true if $i_k=j_k+1$ for $D_r$ or $p_d=d+1$ for $T_r$, so we may assume that $(l+1,l)\notin\Co_l\cap (D\cup T)$ for all $1\leq l\leq r$.

First, consider the case $\Co_{r+1}\cap D=\varnothing$. (And, consequently, $\Co_l\cap D=\varnothing$ for any $k+1\leq l\leq r+1$.) Then $[T]_k\leq[D]_k$ and $t=|[T]_k|=|[D]_k|$. (Note that if $\Co_k\cap T\neq\varnothing$, then $\Co_{k+1}\cap T=\varnothing$.) If $[T]_k=[D]_k$, then $T=D_{(i_k,k)}^{\to}\in N^0(D)$. Suppose $[T]_k<[D]_k$. By the inductive assumption on $n$, there exists $S_1\in N'([D]_k)$ such that $[T]_k\leq S_1$. Then $S=[S_1]^k\in N'(D)$ and $T\leq S$, hence we are done.


iv) Second, consider the case the case $\Co_{r+1}\cap D=\{(i,r+1)\}\neq\varnothing$. If $D\cap\Co_r=\varnothing$ (and so $\Co_l\cap D=\varnothing$ for any $k+1\leq l\leq r$), then we can argue as on the previous step. Assume $D\cap\Co_r\neq\varnothing$, i.e., $k=r$. If $\Co_{r+1}\cap T=\varnothing$, then $D_{r+1}$, $T_{r+1}$ satisfy (\ref{formula:conditions_0}). Suppose $\Co_{r+1}\cap T=\{(p,r+1)\}$. Further, assume $d=|T_r|=r$. Since $T_{r+1}$ does not satisfy (\ref{formula:conditions_0}), $p_r<p$. Put $$\wt T=\left(T\setminus\{(p_r,r),(p,r+1)\}\right)\cup\{(p,r),(p_r,r+1)\}.$$
Obviously, $\wt T>T$ and $|\wt T|=t=|D|$. If $\wt T<D$, then $T\notin\wt L^0(D)$. If $\wt T=D$, then $(i_r,r)=(p,r)\in D$, $(p_r,r+1)=(i,r+1)\in D$ and $S=T=D_{(p,r)}^{(i,r+1)}\in N^0(D)$, so we are done.

But if $\wt T\nleq D$, then $i_r<i$. In this case, we put
$$\wt D=\left(D\setminus\{(i_r,r),(i,r+1)\}\right)\cup\{(i,r),(i_r,r+1)\}.$$
Then $\wt D>D$ and $\wt D>\wt T$, so, by the inductive hypothesis on $\leq$, there exists $S_1\in N'(\wt D)$ such that $\wt T\leq S_1$. On can check that this implies an existence of $S\in N'(D)$ such that $T\leq D$. (In fact, $S$ is obtained from $D$ by the ``same'' operation as $S_1$ from $\wt D$.)

Now, assume $d<r$. Here we put $\wt T=\left(T\setminus\{(p,r+1)\}\right)\cup\{(p,r)\}$. Evidently, $\wt T>T$ and $\wt T\neq D$. If $\wt T<D$, then $T\notin\wt L^0(D)$, so $\wt T\nleq D$. It follows that $i_r<i$. Put
\begin{equation*}
\wt D=\left(D\setminus\{(i_r,r),(i,r+1)\}\right)\cup\{(i,r),(i_r,r+1)\},
\end{equation*}
then $\wt D>D$ and $\wt D>\wt T$. By the inductive hypothesis on $\leq$, there exists $S_2\in N'(\wt D)$ such that $\wt T\leq S_2$. On can check that this implies an existence of $S\in N'(D)$ such that $T\leq D$. (In fact, $S$ is obtained from $D$ by the ``same'' operation as $S_2$ from $\wt D$.) Thus, (\ref{formula:if_conds_then_OK_0}) is proved.

v) From now on, we may assume that $D$, $T$ satisfy (\ref{formula:conditions_0}). First, consider the case $p_d>d+1$. If $i_k=j_k+1$, then $S=D_{(i_k,j_k)}^-\in\wt N^-(D)\subseteq N'(D)$ and $T\leq S$, hence we may assume that $i_k>j_k+1$. If $q_t>j_t$, then, arguing as in the second paragraph of step (iii), one can prove that either $T\in N^0(D)$ or $T\notin\wt L^0(D)$. On the other hand, if $q_t=j_t=t$, then $[D]_t$, $[T]_t$ are rook placements of the same length $t$, and $[D]_t>[T]_t$. The inductive assumption on $n$ shows that there exists $S_1\in N'(D)$ such that $[T]_t\leq S_1$. It is easy to see that $S=[S_1]^t\in N'(D)$ and $T\leq S$, as required.

Second, consider the case $p_d=d+1$. Put $T'=T_{(d+1,d)}^-$.  By the previous Lemma, there exists $S_1\in\wt N^-(D)$ such that $T'\leq S_1$. If $T\leq S_1$, then we are done. If $T\nleq S_1$, then $(i_d,d)\in\wt M(D)$ and $S_1=D'=D_{(i_d,d)}^-$. Denote $T''=[T']_d$, $D''=[D']_d$, then $T''\leq D''$. If $T''=D''$, then $(d+1,d)\notin D$, because $T\neq D$. Hence $D\cap\Ro_{d+1}=\varnothing$, $S=D_{(i_d,d)}^{\uparrow}$ is well-defined and $T\leq S$.

Now, suppose $T''<D''$. By the induction hypothesis on $n$, there exists $S_2\in N'(D'')$ such that ${T''\leq S_2}$. If $S_2\cap\Ro_{i_d}=\varnothing$, then $S=[S_2]^d\cup\{(i_d,d)\}\in N'(D)$ and $T\leq S$. At the contrary, assume $(i_d,j)\in S_2$ for some $j$. If $j<d$, then $S_2=(D'')_{(x,j)}^{\uparrow}$ for some $x$. Since $x>i_d$, we obtain ${(x,j)\in B_{(i_d,d)}(D)}$, $S=D_{(x,j)}^{(i_d,d)}\in N^0(D)$ and $T\leq S$. Finally, if $j>d$, then $S_2=(D'')_{(x,j)}^{\uparrow}$ for some $(x,j)\in D$. In this case, $\wt T=D\setminus\{(x,j)\}>T$ and $|\wt T|=t-1<t=|D|$, so $T\notin\wt L^0(D)$, a~contradiction. The result follows.}

\newpage Finally, we will prove that $L^+(D)=N^+(D)$ (cf. \cite[Lemma 3.7]{Ignatyev1}).

\lemmp{Let $D\in\Ro$. One has $L^+(D)=N^+(D)$.}{By Lemma~\ref{lemm:N_+_in_L_+}, it is enough to check that $L^+(D)\subseteq N^+(D)$. We must show that
\begin{equation}
\begin{split}
&\text{if }
T\leq D\text{ and }|T|>|D|,\\
&\text{then there exists }S\in N''(D)=\wt N^-(D)\cup N^0(D)\cup N^+(D)
\text{ such that }T\leq S<D.
\end{split}\label{formula:proof_+}
\end{equation}
We will proceed by induction on $n$ (for $n=1$, there is nothing to prove). For convenience, we split the proof into five steps.

i) One can easily check that if $D=D_0$, then (\ref{formula:proof_+}) holds. Therefore, we may also use the second (downward) induction on the partial order $\leq$ on $\Ro$.

ii) Let $D<D_0$, $D=\{(i_1,j_1),\ldots,(i_s,j_s)\}$, $T=\{(p_1,q_1),\ldots,(p_t,q_s)\}$, $s<t$ and $T<D$. We claim that
\begin{equation}
\begin{split}
&\text{if (\ref{formula:proof_+}) holds for all $D$, $T$ satisfying (\ref{formula:conditions_0}),}\\
&\text{then (\ref{formula:proof_+}) holds for all $D$, $T\in\Ro$.}
\end{split}
\label{formula:if_conds_then_OK_+}
\end{equation}
Clearly, to prove this, it is enough to show that if (\ref{formula:proof_+}) holds for all $D$, $T$ satisfying (\ref{formula:conditions_0}), and $D_r$,~$T_r$ do not satisfy (\ref{formula:conditions_0}) for some $r$, then (\ref{formula:proof_+}) holds for $D$, $T$. We will proceed by induction on $r$ (the base $r=1$ is clear). Evidently, we may assume that $T\in L^+(D)$.

iii)--iv) The induction step can be performed similarly to steps (iii)--(iv) of the proof of the previous Lemma, so we may assume without loss of generality that $D$, $T$ satisfy (\ref{formula:proof_+}).

v) If $p_d>d+1$, then, arguing as in the first paragraph of step (v) of the proof of the previous Lemma, we obtain the result, hence we may assume $p_d=d+1$. Put $T'=T_{(d+1,d)}^-$. First, consider the case $|D|=s=t-1=|T'|$. By the previous Lemma, there exists $S_1\in N'(D)$ such that $T'\leq S_1$. If $T\leq S_1$, then the result follows. Thus, we may assume that $T\nleq S_1$ (clearly, $T\neq S_1$). If $S_1\in N^0(D)$, then $S_1=D_{(i_d,d)}^{\to}$, so $i_d>d$. Denote $S_1\cap\Ro_{i_d}=\{(i_d,m)\}$, then one can find $\alpha$ such that $d\leq\alpha\leq m$ and $(\alpha,m)\in C_{(i_d,d)}(D)$. Hence $S=D_{(i_d,d)}^{(\alpha,m)}\in N^0(D)$ and $T\leq S$, as required.

Next, suppose $S_1\in\wt N^-(D)$. This means that $(i_d,d)\in\wt M(D)$ and $S_1=D'=D_{(i_d,d)}^-$. Denote $D''=[D']_d$, $T''=[T']_d$, then $D''>T''$ and $|D''|=s-1<|T''|=t-1$. The inductive hypothesis on $n$ says that there exists $S_2\in N''(D'')$ such that $T''\leq S_2$. If $S_2\in N'(D'')$, then one can argue as in the last paragraph of step (v) of the proof of the previous Lemma, so from now on we may assume that $S_2\in N^+(D'')$.

If $S_2=(D'')_{(i,j)}^{i_d,\beta}$ for some $(i,j)\in D''$, then automatically $i>i_d$. If $j>d$, then $S=D_{(i,j)}^{\to}\geq T$ (in fact, $S=\left(D\setminus\{(i,j)\}\right)\cup(i,\beta)$). But if $j\leq d$, then $(i,j-1)\in B_{(i_d,d)}(D)$ and $S=D_{(i_d,d)}^{(i,j-1)}\geq T$. If $S_2=(D'')_{(i,j)}^{\alpha,\beta}$ for some $(i,j)\in D''$ such that $j-1<d$, $(i,j-1)>(i_d,d)$, but $(i_d,d)\nless(\alpha,j-1)$ and $(i_d,d)\nless(i,\beta)$, then $d+1=p_d<p_{j-1}\leq\alpha<i_d$ and $\beta>d$. Hence $\wt T=\left(D\setminus\{(i_d,d)\}\right)\cup\{(\alpha,d)\}>T$ and $T\notin L^+(D)$, a contradiction. In all other cases, $S=[S_2]^d\cup\{(i_d,d)\}\geq T$.

Second, consider the case $|D|=s<t-1=|T'|$. Suppose there exists a root $(\alpha,\beta)\in T$ such that $\alpha\leq i_d$ and $\beta\geq d+1$. Let $\beta$ be the minimal among all such roots, then $D>\wt T>T$, where $$\wt T=\left(T\setminus\{(d+1,d),(\alpha,\beta)\}\right)\cup\{(\alpha,d)\}.$$ Thus, $T\notin L^+(D)$. On the other hand, if such $(\alpha,\beta)$ does not exist, then define $D'$, $T'$, $D''$, $T''$ as above. By the inductive hypothesis on $n$, one can find $S_2\in N''(D'')$ such that $S_2\geq T$ If $S_2\in N'(D'')$, then one can argue as in the last paragraph of step (v) of the proof of the previous Lemma. But if $S_2\in N^+(D'')$, then, arguing as in the previous paragraph, we conclude the proof.}


\newpage\begin{thebibliography}{XXXX}

\bibitem[An1]{Andre1} C.A.M. Andr\`e. Basic sums of coadjoint orbits of the unitriangular group. J.~Algebra \textbf{176} (1995), 959--1000.

\bibitem[An2]{Andre2} C.A.M. Andr\`e. The basic character table of the unitriangular group. J.~Algebra \textbf{241} (2001), 437--471.

\bibitem[AN]{AndreNeto1} C.A.M Andr\`e., A.M. Neto. Super-characters of finite unipotent groups of types $B_n$, $C_n$ and~$D_n$. J. Algebra \textbf{305} (2006), 394--429.

\bibitem[Bu]{Bourbaki} N. Bourbaki. Lie groups and Lie algebras. Chapters 4--6. Springer, 2002.

\bibitem[De]{Deodhar} V. Deodhar. On the root system of a Coxeter group. Comm. Algebra \textbf{10} (1982), no.~6\linebreak 611--630.

\bibitem[Hu1]{Humphreys} J. Humphreys. Linear algebraic groups. Springer, 1975.

\bibitem[Hu2]{Humpreys2} J. Humphreys. Reflection groups and Coxeter groups. Cambridge University Press, Cambridge, 1992.

\bibitem[Ig1]{Ignatyev1} M.V. Ignatyev. Combinatorics of $B$-orbits and the Bruhat--Chevalley order on involutions. Transformation Groups \textbf{17} (2012), no. 3, 747--780, see also arXiv: \texttt{math.RT/1101.2189}.

\bibitem[Ig2]{Ignatyev2} M.V. Ignatyev. The Bruhat--Chevalley order on involutions in the hyperoctahedral group and combinatorics of $B$-orbit closures (in Russian). Zapiski Nauchn. Sem. POMI \textbf{400} (2012), 166-188. English transl.: J. Math. Sci. \textbf{192} (2013), no. 2, 220--231, see also arXiv: \texttt{math.RT/1112.2624}.

\bibitem[Ig3]{Ignatyev3} M.V. Ignatyev. Orthogonal subsets of classical root systems and coadjoint orbits of unipotent groups (in Russian). Mat. Zametki \textbf{86} (2009), no. 1, 65--80. English transl.: Math. Notes \textbf{86} (2009), no. 1, 65--80, see also arXiv: \texttt{math.RT/0904.2841}.

\bibitem[Ig4]{Ignatyev4} M.V. Ignatyev. Orthogonal subsets of root systems and the orbit method (in Russian). Algebra i Analiz 22 (2010), no. 5, 104--130. English transl.: St. Petersburg Math.~J. \textbf{22} (2011), no. 5, 777--794, see also arXiv: \texttt{math.RT/1007.5220}.

\bibitem[In1]{Incitti1} F. Incitti. Bruhat order on the involutions of classical Weyl groups. Ph.D. thesis. Dipartimento di Matematika ``Guido Castelnuovo'', Universit\`a di Roma ``La Sapienza'', 2003.

\bibitem[In2]{Incitti2} F. Incitti. The Bruhat order on the involutions of the symmetric groups. J. Alg. Combin. \textbf{20} (2004), no. 3, 243--261.

\bibitem[Ke]{Kerov} S.V. Kerov. Rook placements on Ferrers boards and matrix integrals. J. Math. Sci. \textbf{96} (1999), no. 5, 3531--3536.

\bibitem[Ki1]{Kirillov1} A.A. Kirillov. Unitary representations of nilpotent Lie groups. Russian Math. Surveys \textbf{17}~(1962), 53--110.

\bibitem[Ki2]{Kirillov2} A.A. Kirillov. Lectures on the orbit method. Grad. Studies in Math. \textbf{64}, AMS, 2004.

\bibitem[Ki3]{Kirillov3} A.A. Kirillov. Two more variations on the
triangular theme. In: C.~Duval, L.~Guieu, V.~Ovsienko (eds.). The
orbit method in geometry and~physics (In honour of A.A. Kirillov).
Progr. in Math. \textbf{213}, Birkh\"auser, Boston, 2003, 243--258.

\bibitem[Ko1]{Kostant1} B. Kostant. The cascade of orthogonal roots and the coadjoint structure of the nilradical of a Borel subgroup of a semisimple Lie group, see arXiv: \texttt{math.RT/1101.5382}.

\bibitem[Ko2]{Kostant2} B. Kostant. Coadjoint structure of Borel subgroups and their nilradicals, see arXiv: \texttt{math.RT/1205.2362}.

\bibitem[Me1]{Melnikov1} A. Melnikov. $B$-orbit in solution to the equation $X^2=0$ in triangular matrices. J. Algebra \textbf{223} (2000), no. 1, 101--108.

\bibitem[Me2]{Melnikov2} A. Melnikov. Description of $B$-orbit closures of order 2 in upper-triangular matrices. Transformation Groups \textbf{11} (2006), no. 2, 217--247.

\bibitem[Pa]{Panov} A.N. Panov. Involutions in $S_n$ and associated coadjoint orbits (in Russian). Zapiski Nauchn. Sem. POMI \textbf{349} (2007), 150--173. English transl.: J. Math. Sci. \textbf{151} (2008), no.~3, 3018--3031, see also arXiv: \texttt{math.RT/0801.3022}.

\bibitem[Pr]{Proctor} R.A. Proctor. Classical Bruhat orders and lexicographical shellability. J. Algebra \textbf{77} (1982), no. 1, 104--126.

\bibitem[RS]{RichardsonSpringer} R.W. Richardson and T.A. Springer, {\it The Bruhat order on symmetric varieties}, Geom. Dedicata \textbf{35} (1990), no. 1--3, 389--436.

\bibitem[Sm]{Smirnov} E. Smirnov. Orbites d'un sous-groupe de Borel dans le produit de deux Grassmanniennes. Ph.D. thesis.
University of Grenoble I, 2007.

\bibitem[Sp]{Springer} T.A. Springer. Some remarks on involutions in Coxeter groups. Comm. Algebra \textbf{10} (1982), no. 6, 631--636.

\end{thebibliography}
\end{document}